%% file: HSIAnomDet_arxiv.tex
\documentclass{article}
\input{myMacros_arxiv.tex}
\usepackage{cite}
\usepackage{verbatim}
\usepackage{url}
\usepackage{subfigure}
\usepackage{amsthm}
\def\PP{{\mathsf P}}

\newtheorem{theorem}{Theorem}
\newtheorem{lemma}{Lemma}
\newtheorem{corollary}{Corollary}

\newcommand{\done}[1]{{#1}}
\usepackage{authblk}

\begin{document}

\title{Target Detection Performance Bounds in Compressive Imaging\thanks{This work was supported by NSF Award No.\ DMS-08-11062, DARPA Grant
   No. HR0011-09-1-0036, and AFRL Grant No. FA8650-07-D-1221.}}



\author[1]{Kalyani Krishnamurthy\thanks{kk63@duke.edu}}
\author[1]{Rebecca Willett\thanks{willett@duke.edu}}
\author[2]{Maxim Raginsky\thanks{maxim@illinois.edu}}
\affil[1]{Department of Electrical and Computer Engineering, Duke University, Durham, NC 27708.}
\affil[2]{Department of Electrical and Computer Engineering and Coordinated Science Laboratory, University of Illinois at Urbana-Champaign, Urbana, IL 61801.}

\renewcommand\Authands{ and }

\date{\today}

\maketitle

\begin{abstract}
  This paper describes computationally efficient approaches and
  associated theoretical performance guarantees for the detection of
  known targets and anomalies from
  few projection measurements of the underlying signals. The proposed approaches accommodate signals of different strengths contaminated by a colored Gaussian background, and perform detection without reconstructing the underlying signals from the observations.
The theoretical performance bounds of the target
  detector highlight fundamental tradeoffs among the number of
  measurements collected, amount of background signal present, signal-to-noise ratio, and 
  similarity among potential targets coming from a known dictionary. The anomaly
  detector is designed to control the number of false discoveries.
  The proposed approach does not depend on a known sparse
  representation of targets; rather, the theoretical performance
  bounds exploit the structure of a known dictionary of targets and {the distance preservation property of the measurement matrix.}
  Simulation experiments illustrate
  the practicality and effectiveness of the proposed approaches.
  \end{abstract}


\sloppypar
\section{Introduction}
\label{sec:intro}
{The theory of compressive sensing (CS) has shown that it is possible
to accurately \emph{reconstruct} a sparse signal from few (relative to the signal dimension)
projection measurements \cite{CS:candes2,CS:donoho}. Though such a reconstruction is crucial to visually inspect the signal, 
there are many instances where one is solely interested in identifying whether the underlying signal is one of several possible signals of interest. In such situations, a complete reconstruction is computationally expensive and does not optimize the correct performance metric. Recently, CS ideas have been exploited in \cite{Baraniuk_Smashed_filter,Nowak_CS_ForSignalClassification,Baraniuk_SparseSignalDet} to perform target detection and classification from projection measurements, without reconstructing the underlying signal of interest.  In \cite{Baraniuk_Smashed_filter,Nowak_CS_ForSignalClassification}, the authors propose nearest-neighbor based methods to classify a signal $\f \in \reals^N$ to one of $m$ known signals given projection measurements of the form $\by = \bA\f + \bn \in \reals^K$ for $K \leq N$, where $\bA \in \reals^{K \times N}$ is a known projection operator and $\bn \sim \Gaussian{\zeros,\sigma^2\bI}$ is the additive Gaussian noise. This model is simple to analyze, but is impractical, since in reality, a signal is always corrupted by some kind of interference or background noise. Extension of the methods in \cite{Baraniuk_Smashed_filter,Nowak_CS_ForSignalClassification} to handle background noise is nontrivial. Though \cite{Baraniuk_SparseSignalDet} provides a way to account for background contamination, it makes a strong assumption that the signal of interest and the background are sparse in bases that are incoherent. This might not always be true in many applications. 
Recent works on CS \cite{aeron2010information,arias2011noise} allow for the input signal $\f$ to be corrupted by some pre-measurement noise $\bb \sim \Gaussian{\zeros,\sigma_b^2\bI}$ such that one observes $\by = \bA(\f+\bb) + \bn$, and study reconstruction performance as a function of the number of measurements, pre- and post-measurement noise statistics and the dimension of the input signal. In this work, however, we are interested in performing target detection without an intermediate reconstruction step. Furthermore, the increased utility of high-dimensional imaging techniques such as spectral imaging or videography in applications like remote sensing, biomedical imaging and astronomical
imaging \cite{forest_remsensing,stellman2000real,rhessi,johnson:014036,martin1998determining,martin2006development,han2007fusion,zuzak2008intraoperative} necessitates the extension of compressive target detection ideas to such imaging modalities to achieve reliable target detection from fewer measurements relative to the ambient signal dimensions. }

{For example, recent advances in compressive sensing (CS) have led to the development of
new spectral imaging platforms which attempt to address challenges in conventional imaging platforms 
related to system size, resolution, and noise by acquiring fewer compressive
measurements than spatiospectral voxels \cite{Wagadarikar:08,hyperspectralCS_BradyBecca,brady:62460A,riceCamera,woolfe2006hyper,SLM}. However, these system designs have a number of
degrees of freedom which influence subsequent data analysis.
For instance, the single-shot compressive spectral imager discussed in \cite{hyperspectralCS_BradyBecca} 
collects \emph{one} coded projection of each spectrum in the scene. One projection per spectrum is sufficient for reconstructing spatially homogeneous spectral images, since projections of neighboring locations can be combined to infer each spectrum. Significantly more projections are required for detecting targets of unknown strengths without the benefit of spatial homogeneity. We are interested in investigating how several such systems can be used in \emph{parallel} to reliably detect spectral targets and anomalies from different coded projections.}

In general, we consider a broadly applicable framework that allows us to account for background and sensor noise, and perform target detection directly from projection measurements of signals obtained at different spatial or temporal locations. The precise problem formulation is provided below.

\subsection{Problem formulation}
\label{subsec:problem formulation}
Let us assume access to a dictionary of possible targets of interest $\sD =
\{\f^{(1)},\f^{(2)},\ldots,\f^{(m)}\}$, where $\f^{(j)} \in \reals^N$
for $j = 1, \ldots, m$ is unit-norm. Our measurements are of the form
\begin{align}
\bz_i = \bPhi (\alpha_i \f^{*}_i+\bb_i) + \bw_i
\label{eqn:problem setup}
\end{align}
where 
\begin{itemize}
\item $i \in \{1,\ldots,M\}$ indexes the spatial or temporal locations at which data are collected;
\item $\alpha_i \geq 0$ is a measure of the signal-to-noise ratio
  at location $i$, which is either known or estimated from observations;
\item $\bPhi \in \reals^{K \times N}$ for $K < N$, is a measurement matrix to be specified in Sec.~\ref{sec:preprocessing};
\item $\bb_i \in \reals^{N} \sim \sN(\bmu_b,\bSigma_b)$ is the background noise vector, and 
$\bw_i \in \reals^K \sim \sN(\zeros,\sigma^2 \bI)$ is the
  i.i.d.\ sensor noise.
\end{itemize}
For example, in the case of spectral imaging $\fstar_i$ represents the spectrum at the \ith spatial location, and in video sequences $\fstar_i$ represents the vectorized image frame obtained at the \ith time interval. 
In this paper we consider the following target detection
problems:
\begin{enumerate}
\item Dictionary signal detection (DSD): Here we assume that each
$\fstar_i \in \sD$ for $i \in \setbr{1,\ldots,M}$, and our task
is to detect all instances of one target signal $\f^{(j)} \in \sD$ for some unknown $j \in \setbr{1,\ldots,m}$, \ie to locate $S = \setbr{i: \f^{*}_i = \f^{(j)}}$. DSD is useful in contexts in which we know the makeup of a scene and wish to focus our attention on the locations of
a particular signal. For instance, in spectral imaging, DSD is used to study a scene of interest by classifying
every spectrum in the scene to different known classes \cite{martin1998determining,manolakis2003hyperspectral}. In a video setup, DSD could be used to classify video segments to one of several categories (such as news, weather, sports, etc.) by projecting the video sequence to an appropriate feature space and comparing the feature vectors to the ones in a known dictionary \cite{tvVideoClassification}. 
\item Anomalous signal detection (ASD): Here, our task is to detect all
signals which are {\em not} members of our dictionary, \ie detect $S = \setbr{i:\f^{*}_i \notin \sD}$. (This is akin to anomaly detection methods in the literature which are based on nominal, nonanomalous training samples \cite{Hero06geometricentropy,Steinwart05aclassification}.)  
For instance, ASD may be used when we know most components of a spectral image and wish to identify all spectra which
deviate from this model \cite{stein2002adh}.
\end{enumerate}  

{Our goal is to accurately perform DSD or ASD without reconstructing
the spectral input $\f^{*}_i$ from $\bz_i$ for $i \in \{1,\ldots,M\}$. Accounting for background is a crucial issue. Typically, the background corresponding to the scene of interest and
the sensor noise are modeled together by a colored multivariate
Gaussian distribution \cite{manolakis2002dah}. However, in our case,
it is important to distinguish the two because of the presence of the
projection operator $\bPhi$. The projection operator acts upon the
background spectrum in the same way as on the target spectrum, but it
does not affect the sensor noise. We assume that $\bb_i$ and
$\bw_i$ are independent of each other, and the prior probabilities of
different targets in the dictionary $p^{(j)} = \prob\left(\f^{*}_{i}=\f^{(j)}\right)$ for $j \in \{1,\cdots,m\}$ are known in
advance. If these probabilities are
unknown, then the targets can be considered equally likely. Given this
setup, our goal is to develop suitable target and anomaly detection
approaches, and provide theoretical guarantees on their performances.
}

In this paper we develop detection performance bounds which show how performance
scales with the number of detectors in a compressive setting as a
function of SNR, the similarity between potential
targets in a known dictionary, and their prior probabilities.  
Our bounds are based on a detection strategy which operates directly
on the collected data as opposed to first reconstructing each $\fstar_i$ and then performing detection on the estimated signals. Reconstruction as an intermediate step in detection may be
appealing to end users who wish to visually inspect spectral images
instead of relying entirely on an automatic detection
algorithm. However, using this intermediate step has two potential
pitfalls. First, the Rao--Blackwell theorem \cite{Berger} tells
us that an optimal detection algorithm operating on the processed
data (\ie not sufficient statistics) {\em cannot} perform better than an optimal detection algorithm
operating on the raw data. In other words, optimal performance is
possible on the raw data, but we have no such performance guarantee
for the reconstructed signals. Second, the relationship between reconstruction errors and detection performance is not well understood in many settings. Although we do not reconstruct the underlying signals, our performance bounds are intimately related to the signal resolution needed to achieve the signal diversity present in our dictionary. Since we have many fewer observations than the signals at this resolution, we adopt the ``compressive" terminology.

\subsection{Performance metric}
To assess the performance of our detection strategies, we consider the False Discovery Rate (FDR) metric and related quantities developed for multiple hypothesis testing problems \cite{benjamini1995cfd}. Since we collect $M$ independent observations of potentially different signals, we are simultaneously conducting $M$ hypothesis tests when we search for targets. Unlike the probability of false alarm, which measures the probability of falsely declaring a target for a \emph{single} test, the $\FDR$ measures the \emph{fraction} of declared targets that are false alarms, that is, it provides information about the entire set of $M$ hypotheses instead of just one. More formally, the $\FDR$ is given by, 
\begin{align*}
  \FDR = \expect\left[\frac{V}{R}\right],
\end{align*}
where $V$ is the number of falsely rejected null hypotheses, and
$R$ is the total number of rejected null hypotheses. Controlling the false discovery rate in a multiple
hypothesis testing framework is akin to designing a constant false
alarm rate (CFAR) detector in spectral target detection applications
that keeps the false alarm rate at a desired level irrespective of the
background interference and sensor noise statistics
\cite{manolakis2003hyperspectral}. 

\subsection{Previous investigations} 
Much of the classical target detection literature \cite{kelly1986ada,adaptiveSubsDet,scharf1994matched,jin2009comparative,kwon2006kernel} assume that each target lies in a $P$-dimensional subspace of $\reals^N$ for $P < N$. 
The subspace in
which the target lies is often assumed to be known or specified by the user, and the variability of the background is modeled
using a probability distribution. Given knowledge of the target
subspace, background statistics and sensor
noise statistics, detection methods based on LRTs (likelihood ratio tests) and GLRTs (generalized likelihood ratio tests) have been proposed in 
\cite{kelly1986ada,scharf1996adaptive,adaptiveSubsDet,scharf1994matched,jin2009comparative,kwon2006kernel}. A subspace model is optimal if the subspace in which targets lie is known in advance. However, in many applications, such subspaces might be hard to characterize. An alternative, and a more flexible option is to assume that the high-dimensional target exhibits some low-dimensional structure that can be exploited to perform efficient target detection. This approach is utilized in this work and in \cite{Nowak_CS_ForSignalClassification} where the target signal in $\reals^N$ is assumed to come from a dictionary of $m$ known signals such that $m \ll N$, and in \cite{Baraniuk_Smashed_filter}, where the targets are assumed to lie in a low-dimensional manifold embedded in high-dimensional target space.

Recently, several methods for target or anomaly detection that rely on recovering the full
spatiospectral data from projection measurements \cite{parmar2008spatio,WillettGehm07} have been proposed. However, they are computationally intensive and the detection performance associated with these reconstructions is unknown. 
Other researchers have
exploited compressive sensing to perform target detection and
classification without reconstructing the underlying signal
\cite{Baraniuk_Smashed_filter,Nowak_CS_ForSignalClassification,Baraniuk_SparseSignalDet}.
In \cite{Baraniuk_SparseSignalDet}, the authors propose a matching
pursuit based algorithm, called the \emph{Incoherent Detection and
  Estimation Algorithm} (IDEA), to detect the presence of a signal of
interest against a strong interfering signal from noisy projection measurements. The algorithm is shown to perform well on experimental
data sets under some strong assumptions on the sparsity of the signal
of interest and the interfering signal. 
In
\cite{Baraniuk_Smashed_filter}, the authors develop a classification
algorithm called the \emph{smashed filter} to classify an image in
$\reals^{N}$ to one of $m$ known classes from $K$ 
projections of the signal, where $K < N$. The underlying image is
assumed to lie on a low-dimensional manifold, and the algorithm finds
the closest match from the $m$ known classes by performing a nearest
neighbor search over the $m$ different manifolds. The projection measurements are chosen to preserve the distances among the
manifolds. Though \cite{Baraniuk_Smashed_filter} offers theoretical
bounds on the number of measurements necessary to preserve distances among different manifolds, it is not clear how the performance scales with $K$ or how to incorporate background models into this setup. Moreover, this
approach may be computationally intensive since it involves 
learning and searching over different manifolds. 
In \cite{Nowak_CS_ForSignalClassification}, the authors use a
nearest-neighbor classifier to classify an $N$-dimensional signal to
one of $m$ equally likely target classes based on $K < N$ random
projections, and provide theoretical guarantees on the detector
performance. While the method discussed in \cite{Nowak_CS_ForSignalClassification} is
computationally efficient, it is nontrivial to extend to the
case of target detection with colored background noise
and nonequiprobable targets. Furthermore, their performance guarantees cannot be directly extended to our problem since we focus on error measures that let us analyze the performance of multiple hypothesis tests \emph{simultaneously} as opposed to the above methods that consider compressive classification performance for a {\em single} hypothesis test.  

\done{The authors of a more recent work \cite{fowler2011anomaly} extend the classical RX anomaly detector \cite{RXAlgo} to \emph{directly} detect anomalies from random, orthonormal projection measurements without an intermediate reconstruction step. They numerically show how the detection probability improves as a function of the signal-to-noise ratio when the number of measurements changes. Though probability of detection is a good performance measure, in many applications controlling the false discoveries below a desired level is more crucial. As a result, in our work, we propose an anomaly detection method that controls the false discovery rate below a desired level.}



{
\subsection{Contributions}
This paper makes the following contributions to the above literature:
\begin{itemize}
\item A \emph{compressive target detection} approach, which (a) is computationally efficient, (b) allows for the signal strengths of
  the targets to vary with spatial location, (c) allows for
  backgrounds mixed with potential targets, (d) considers
  targets with different a priori probabilities, and (e) yields
  theoretical guarantees on detector performance.
  This paper unifies preliminary work by the authors \cite{icassp2010TargetDet,icip2010TargetDet}, presents previously unpublished aspects of the proofs, and contains updated experimental results. 
  \item A computationally efficient \emph{anomaly detection} method that
  detects  anomalies of different strengths from projection measurements and also controls the false discovery rate at a
  desired level.
  \item A \emph{whitening filter approach} to compressive measurements of signals with background contamination, and associated analysis leading to bounds on the amount of background to which our detection procedure is robust.
\end{itemize}
The above theoretical results, which are the main focus of this paper, are supported with simulation studies in Sec.~\ref{Sec:ExperimentalResults}.
Classical detection methods described in \cite{manolakis2002dah,manolakis2003hyperspectral,boardman1993spectral,scharf1996adaptive,kelly1986ada,adaptiveSubsDet,scharf1994matched,guoa-template,szlam-split,jin2009comparative,kwon2006kernel,stein2002adh,
  RXAlgo,KernelRXAnomDet} do not establish performance bounds as a function of signal resolution or target dictionary properties and rely on relatively \emph{direct observation} models which we show to be suboptimal when the detector size is limited. The methods in \cite{Baraniuk_Smashed_filter} and \cite{Baraniuk_SparseSignalDet} do not contain performance analysis, and our analysis builds upon the analysis in \cite{Nowak_CS_ForSignalClassification} to account for several specific aspects of the compressive target detection problem.
}

\section{Whitening compressive observations}
\label{sec:preprocessing}
Before we present our detection methods for DSD and ASD problems respectively, we briefly discuss a whitening step that is common to both our problems of interest. 

{Let us suppose that there are enough background
training data available to estimate the background mean $\bmu_{b}$ and
covariance matrix $\bSigma_b$. We can assume without loss of
generality that $\bmu_{\bb} = \zeros$ since $\bPhi \bmu_{\bb}$ can be
subtracted from $\by$. Given the knowledge of the background statistics, we can transform the background and sensor noise model $\bPhi \bb_i + \bw_i\sim \sN(\zeros, \bPhi \bSigma_b \bPhi^T + \sigma^2\bI)$ discussed in \eref{eqn:problem setup} to a simple white Gaussian noise model by multiplying the observations $\bz_i$, $i \in \{1,\ldots,M\}$, by the \emph{whitening filter} $C_{\bPhi} \triangleq (\bPhi \bSigma_b \bPhi^T + \sigma^2\bI)^{-1/2}.$
This whitening transformation reduces the observation model in \eref{eqn:problem setup} to   
\begin{align}
&\by_i =  C_{\bPhi}\underbrace{\left(\bPhi\left(\alpha_i\f^{*}_{i} + \bb_i\right) + \bw_i\right)}_{\bz_i}= \alpha_i \bA \f^{*}_i + \bn_i \label{eqn:whitened}
\end{align}
where 
\begin{align}
\bA = C_{\bPhi}\bPhi, \label{eqn:CPhiPhi}
\end{align}
and $\bn_i = C_{\bPhi}\left(\bPhi \bb_i + \bw_i\right) \sim \sN(\zeros,\bI)$. To verify that $\bn_i \sim \sN(\zeros,\bI)$, observe that $$\bn_{i} = C_{\bPhi}\left(\bPhi \bb_i + \bw_i\right) \sim \sN\Big(\zeros,\underbrace{C_{\bPhi}\left(\bPhi \bSigma_b \bPhi^T + \sigma^2\bI\right)C_{\bPhi}^{T}}_{\bI}\Big).$$ 
We can now choose $\bPhi$ so that the corresponding $\bA$ has certain desirable properties as detailed in Sec.~\ref{sec:targetDetection} and Sec.~\ref{sec:AnomDet}.}
{For a given $\bA$, the following theorem provides a construction of $\bPhi$ that satisfies \eref{eqn:CPhiPhi} and a bound on the maximum tolerable background contamination:}
\\
{
\begin{theorem} 
Let $\bB = \bI - \bA \bSigma_b \bA^T$. If the largest eigenvalue of $\bSigma_b$ satisfies 
\begin{align}
\lambda_{\max} < \frac{1}{\|\bA\|^2}, \label{eqn:bgdCond}
\end{align}
where $\|\bA\|$ is the spectral norm of $\bA$, then $\bB$ is positive definite and $\bPhi = \sigma \bB^{-1/2}\bA$ is a sensing matrix, which can be used in conjunction with a whitening filter to produce observations modeled in \eref{eqn:whitened}.
\label{thm:PhiConstruction}
\end{theorem}
}
The proof of this theorem is provided in Appendix~\ref{app:PhiConstruction}. This theorem draws an interesting relationship between the maximum background perturbation that the system can tolerate and the spectral norm of the measurement matrix, which in turn varies with $K$ and $N$. Hardware designs such as those in \cite{SLM,riceCamera} use spatial light modulators and digital micro mirrors, which allow the measurement matrix $\bPhi$ to be adjusted easily in response to changing background statistics and other operating conditions.

{In the sections that follow, we consider collecting measurements of the form $\by_i = \alpha_i \bA
\f^{*}_i + \bn_i$ given in \eref{eqn:whitened}, where $\f^{*}_i$ is
the target of interest for $i = 1,\ldots,M$, and $\bA \in \reals^{K \times N}$ is a sensing matrix that satisfies \eref{eqn:CPhiPhi}. It is assumed that any background contamination
has been eliminated with the whitening procedure described in this section.}

\section{Dictionary signal detection}
\label{sec:targetDetection}
Suppose that the end user wants to test for the presence of one known
target versus the rest, but it is not known a priori which target from
$\sD$ the user wants to detect. In this case, let us cast the
DSD problem as a multiple hypothesis
testing problem of the form
\begin{align}
\sH^{(j)}_{0i}: \f^{*}_i = \f^{(j)} \quad \text{ vs. } \quad \sH^{(j)}_{1i}: \f^{*}_i \neq \f^{(j)} \label{eqn:Hypothesis test H1}
\end{align}
where $\f^{(j)} \in \sD$ is the target of interest and $i = 1,\ldots,M$. 

\subsection{Decision rule}
\label{subsec:DSD_DecisionRule}
We define our decision rule corresponding to target $\f^{(j)} \in \sD$ in terms of a set of significance regions $\Gamma_{i}^{(j)}$ such that one rejects the $\ith$ null hypothesis if its test statistic $\by_{i}$ falls in the $\ith$ significance region. Specifically, $\Gamma_{i}^{(j)}$ is defined according to
\begin{align}
&\Gamma_{i}^{(j)} = \Big\{\by: \log\prob\big( \f_i^{*}=\f^{(j)} \big|\by_i,\alpha_{i},\bA\big) \leq \label{eqn:ourSigReg}\\
&\log\prob\big(\f_i^{*}=\f^{(\ell)}\big|\by_i,\alpha_{i},\bA\big) \text{ for some } \ell \in \{1,\ldots,m\}, \l \neq j \Big\},\nonumber
\end{align}
where 
$\log \prob\big(\f_i^{*}=\f^{(j)}\big|\by_i,\alpha_{i},\bA\big) = \frac{K}{2}\log\left(\frac{1}{2\pi}\right) -
\frac{\left\|\by_i-\alpha_i\bA\f^{(j)}\right\|^2}{2} + \log p^{(j)}$
is the logarithm of the a posteriori probability density of the target $\f^{(j)}$ at the $i\th$ spatial location given the observations $\by_i$, the signal-to-noise ratio $\alpha_{i}$ and the sensing matrix $\bA$, and $p^{(j)}$ is the a priori probability of target class $j$. \done{Note that the process of determining these decision regions involves a sequence of nearest-neighbor calculations, so the computational complexity scales with the number of classes $m$. In this work, we operate under the assumption that $m$ is much smaller than the dimensionality of the datasets we consider. For example, if we consider spectral images, then the number of objects (signal classes) that make up a scene of interest is often smaller than the number of voxels in the image. This assumption is not unrealistic and has been exploited in earlier work such as \cite{manolakis2003hyperspectral} and the references therein. In most of the prior work we have surveyed \cite{chang10virtual,estSigClass_HSI}, the number of signal classes is less than 35, which doesn't make our approach intractable. }

The decision rule can be formally expressed in terms of the significance regions as follows:
\begin{align}
\text{reject } \sH_{0i}^{(j)} \text{ if the test statistic } \by_{i} \in \Gamma_{i}^{(j)}. \label{eqn:DecisionRule}
\end{align}

We analyze this detector by extending the positive False Discovery Rate ($\pFDR$) error measure introduced by Storey to characterize the errors encountered in performing multiple, independent and \emph{nonidentical} hypothesis tests simultaneously \cite{storey2003positive}. The $\pFDR$, discussed formally below, is the fraction of falsely rejected null hypotheses among the total number of rejected null hypotheses, subject to the positivity condition that one rejects at least one null hypothesis. The $\pFDR$ is similar to the $\FDR$ except that the positivity condition is enforced here. In our context, the positivity condition means that we declare at least one signal to be a nontarget, which in turn implies that the scene of interest is composed of more than one object in the case of spectral imaging, or that the scene is not static in the case of video imaging.

Consider a collection
of significance regions $\bm{\Gamma} = \big\{\Gamma_i^{(j)} : i = 1,\cdots,M \big\}$, such that one declares $\sH_{1i}^{(j)}$ if the test
statistic $\by_i \in \Gamma_i^{(j)}$. The pFDR for multiple,
\emph{nonidentical} hypothesis tests can be defined in terms of the
significance regions as follows:
\begin{align}
\pFDR^{(j)}\left(\bm{\Gamma}\right) &= \expect \left[\left.\frac{V\left(\bm{\Gamma}\right)}{R\left(\bm{\Gamma}\right)} \right| R\left(\bm{\Gamma}\right) > 0\right] \label{eqn:pFDRdefn}
\end{align}
where
\begin{align}
V\left(\bm{\Gamma}\right) &= \sum_{i=1}^{M} \ind{\by_i \in \Gamma_i^{(j)}} \ind{\sH_{0i}} \label{eqn:Vdefn}
\end{align}
is the number of falsely rejected null hypotheses, 
\begin{align}
R\left(\bm{\Gamma}\right) &= \sum_{i=1}^{M} \ind{\by_i \in \Gamma_i^{(j)}}
\end{align}
is the total number of rejected null hypotheses, and $\ind{E} = 1$ if event $E$ is true and $0$ otherwise.
In our setup, the pFDR corresponds to the expected ratio of the number of missed targets to the number of signals declared to be nontargets subject to the condition that at least one signal is declared to be a nontarget. (Note that this ratio is traditionally referred to as the positive false \emph{nondiscovery} rate ({\rm{pFNR}}), but is technically the {\rm{pFDR}} in this context because of our definitions of the null and alternate hypotheses.) 
The theorem below presents our main result:
{
\begin{theorem}
\label{thm:pFDR}
Given observations of the form \eref{eqn:whitened}, if one performs multiple, independent, nonidentical hypothesis tests of the form \eref{eqn:Hypothesis test H1} and decides according to \eref{eqn:DecisionRule}, then the worst-case $\pFDR$ given by $\pFDR_{\max} = \max_{j \in \{1,\ldots,m\}} \pFDR^{(j)}\left(\bm{\Gamma}\right),$ satisfies the following bound:
\begin{align}
\pFDR_{\max}\leq \min\left(1,\frac{\pe_{\max}}{1-\pmax-\pe_{\max}}\right) \label{eqn:pfdr_max_anyA}
\end{align}
where 
\begin{align}
\pmax &= \max_{j \in \{1,\ldots,m\}}p^{(j)}, \nonumber\\
\pe_{\max} &= \max_{i \in \{1,\ldots,M\}} \prob\left(\bhf_i \neq \f^{*}_i\right), \text{ and }\nonumber \\
\bhf_i &= \argmax_{\f \in \sD}\prob\left(\left.\f_i^{*}=\f\right|\by_i,\alpha_{i},\bA\right). \label{eqn:fhat}
\end{align}
\end{theorem}
The proof of this theorem is detailed in Appendix~\ref{app:pFDR}. A key element of our proof is the adaptation of the techniques from \cite{storey2003positive} to {\em nonidentical} independent hypothesis tests.}

\subsection{An achievable bound on the worst-case \pFDR}

Theorem~\ref{thm:pFDR} in the preceding section shows that, for a given $\bA$, the worst-case \pFDR is bounded from above by a function of the worst-case misclassification probability. In this section, we use this theorem to establish an achievable bound on the worst-case \pFDR that explicitly depends on the number of observations $K$, signal strengths $\{\alpha_i\}_{i=1}^{M}$, similarity among different targets of interest, and a priori target probabilities.

Let us first define the quantities
\begin{align*}
	\dmin &= \min_{\f^{(i)},\f^{(j)} \in \sD, i\neq j}\|\f^{(i)}-\f^{(j)}\|\\
	\pmin &= \min_{j \in \setbr{1,\ldots,m}} p^{(j)}\\
	\alpha_{\min} &= \min_{i \in \{1,\ldots,M\}} \alpha_i.
\end{align*}
Then we have the following theorem, whose proof is given in Appendix~\ref{app:achievability}:


\begin{theorem}\label{thm:achievability} Let $\lambda_{\max}$ denote the largest \done{eigenvalue} of $\bSigma_b$.
For a given $0 < \epsilon < 1-\pmax$, assume that $K$ and $N$ are sufficiently large so that the following conditions hold:
\begin{subequations}
\begin{align}
 1 - \pmax - \epsilon &\ge \frac{1-\pmin}{\pmin} \left(1+\frac{\done{\alpha_{\min}}^2 \dmin^2}{4K\sigma^2}\right)^{-\frac{K}{2}} + 2\exp\left(-\frac{(K+N)\epsilon^2}{2}\right)   \label{eq:positive_prob}\\
	\lambda_{\max} &< \frac{1}{(1+\epsilon)^2\left(\sqrt{\frac{N}{K}}+1\right)^2},
		\label{eq:weak_background}\\
	K &> \frac{2\log\left(\frac{2}{\pmin} \frac{1-\pmin}{1-\pmax}\right) }{\log \left(1+\frac{\alpha_{\min}^2 d^2_{\min}}{4K}\right)}. \label{eqn:boundonK}
	\end{align}
	\end{subequations}
Then there exists a $K \times N$ sensing matrix $\bA$ that satisfies the condition of Theorem~\ref{thm:PhiConstruction}, and for which
\begin{align}
	&\pFDR_{\max} \le \frac{1}{\pmin}\left(\frac{1-\pmax}{1-\pmin}\left(1+\frac{\alpha_{\min}^2 \dmin^2}{4K}\right)^{\frac{K}{2}}-\frac{1}{\pmin}\right)^{-1}  + \nonumber \\
	&\qquad \frac{2(1-\pmax)}{\epsilon^2}\exp\left(-\frac{(K+N)\epsilon^2}{2}\right). \label{eq:achievable_pfdr}
\end{align}
\end{theorem}
This result has the following implications and consequences:
\begin{enumerate}
	\item For a given $N$, the upper bound \eqref{eq:weak_background} on $\lambda_{\max}$ increases as $K$ increases, which implies that the system can tolerate more background perturbation if we collect more measurements. 
\item The \pFDR bound \eqref{eq:achievable_pfdr} decays with the increase in the values of $K$, $\dmin$ and $\alpha_{\min}$, and increases as $\pmin$ decreases. For a fixed $\pmax$, $\pmin$, $\alpha_{\min}$ and $\dmin$, the bound in \eref{eq:achievable_pfdr} enables one to choose a value of $K$ to guarantee a desired \pFDR value.
\item The dominant part of the bound \eqref{eq:achievable_pfdr} is independent of $N$, and is only a function of $K$, $\pmax$, $\pmin$, $\alpha_{\min}$, and $\dmin$. The lack of dependence on $N$ is not unexpected. Indeed, when we are interested in preserving pairwise distances among the members of a fixed dictionary of size $m$, the Johnson--Lindenstrauss lemma \cite{originalJL} says that, with high probability, $K = \Order\left(\log m\right)$ random Gaussian projections suffice, regardless of the ambient dimension $N$. This is precisely the regime we are working with here.
\item The bound on $K$ given in \eref{eqn:boundonK} increases logarithmically with the increase in the difference between $\pmax$ and $\pmin$. This is to be expected since one would need more measurements to detect a less probable target as our decision rule weights each target by its a priori probability. If all targets are equally likely, then $\pmax = \pmin = 1/m$, and $K = \Order\left(\log m\right)$ is sufficient provided $\alpha_{\min}^2\dmin^2$ is sufficiently large such that
\begin{align*}
	\log \left(1+\frac{\alpha_{\min}^2 d^2_{\min}}{4K}\right) > \log \left(1+\frac{\alpha_{\min}^2 d^2_{\min}}{4N}\right) >1
	\end{align*}
(where the first inequality holds since $K < N$). In addition, the lower bound on $K$ also illustrates the interplay between the signal strength of the targets, the similarity among different targets in $\sD$, and the number of measurements collected. A small value of $\dmin$ suggests that the targets in $\sD$ are very similar to each other, and thus $\alpha_{\min}$ and $K$ need to be high enough so that similar targets can still be distinguished. The experimental results discussed in Sec.~\ref{Sec:ExperimentalResults} illustrate the tightness of the theoretical results discussed here.
\end{enumerate}

Inspection of the proof shows that if $\bA$ is generated according to a Gaussian distribution, then the conditions of Theorem~\ref{thm:achievability} will be met with high probability.

\section{Extension to a manifold-based target detection framework}
The DSD problem formulation in Sec.~\ref{subsec:problem formulation} is accurate if the signals in the dictionary are faithful representations of the target signals that we observe. In reality, however, the target signals will differ from the dictionary signals owing to the differences in the experimental conditions under which they are collected. For instance, in spectral imaging applications, the observed spectrum of any material will not match the reference spectrum of the same material observed in a laboratory because of the differences in atmospheric and illumination conditions. 
To overcome this problem, one could form a large dictionary to account for such uncertainties in the target signals and perform target detection according to the approaches discussed in Sec.~\ref{sec:preprocessing} and Sec.~\ref{sec:targetDetection}. A potential drawback with this approach is that our theoretical performance bound increases with the size of $\sD$ through $\pmin$ and $\dmin$. Instead, one could reasonably model the target signals observed under different experimental conditions to lie in a low-dimensional submanifold of the high-dimensional ambient signal space as shown to be true for spectral images in \cite{healey1999models}. 
We can exploit this result to extend our analysis to a much broader framework that accounts for uncertainties in our dictionary.

Let us consider a dictionary of manifolds $\sD_{\sM} = \left\{\sM^{(1)},\ldots,\sM^{(m)}\right\}$ corresponding to $m$ different target classes, and that $\fstar_i$ for $i \in \setbr{1,\ldots,M}$ is in one of the manifolds in $\sD_{\sM}$.
Considering an observation model of the form given in \eref{eqn:whitened}, our goal is to determine $\setbr{i: \f^{*}_i \in \sM^{(j)}}$, where $j \in \setbr{1,\ldots,m}$ is the target class of interest. Let us assume that all target classes are equally likely to keep the presentation simple, though the analysis extends to the case where the targets classes have different a priori probabilities. Suppose that we collect independent sets of measurements $\setbr{\by_i}_{i=1}^M$ and $\setbr{\bty_i}_{i=1}^M$. Then, we can use the following two-step procedure to extend our DSD method to this manifold-based framework:
\begin{enumerate}
\item Given $\setbr{\by_i}$, form a data-dependent dictionary $\sD_{\by_i} = \left\{\btf_i^{(1)},\ldots,\btf_i^{(m)}\right\}$ corresponding to each $\by_i$ by finding its nearest-neighbor in each manifold:
$$\btf^{(\ell)}_i = \argmax_{\f \in \sM^{(\ell)}} \prob\left(\left.\by_i\right|\f^{*}_i = \f,\alpha_i,\bA\right)$$ for $\ell \in \{1,\ldots,m\}$ and $i = 1,\ldots,M$.\\
\item Given $\setbr{\bty_i}$ and corresponding $\setbr{\sD_{\by_i}}$, find $$\bhf_i = \argmax_{\btf \in \sD_{\by_i}} \prob\left(\left.\bty_i\right|\f^{*}_i = \btf,\alpha_i,\bA\right)$$ and 
declare that the $\ith$ observed spectrum corresponds to class $j$ if $\bhf_i = \btf_i^{(j)}$.
\end{enumerate}
This two-step procedure is studied in \cite{Baraniuk_Smashed_filter} for the case $\setbr{\by_i} = \setbr{\bty_i}$ where the authors provide bounds on the number of projection measurements needed to preserve distances among manifolds. However, they do not offer associated target detection performance guarantees. Our analysis and the theoretical performance bounds extend directly to this framework if we collect two sets of observations as discussed above. Specifically, the hypothesis tests corresponding to the second step can be written as $$\sH_{0i}: \f^{*}_i = \btf_i^{(j)} \text{ vs. } \sH_{1i}: \f^{*}_i \neq \btf_i^{(j)}$$ where $\btf_i^{(j)} \in \sD_{\by_i}$ for $i = 1,\ldots,M$. Since the dictionary in this case changes with $i$, these tests are \emph{nonidentical}. This is another instance where our extension of pFDR-based analysis towards simultaneous testing of multiple, independent, and nonidentical hypothesis tests \eref{eqn:pFDRdefn} is very significant. Following the proof techniques discussed in the appendix, we can straightforwardly show that  the bound in \eref{eq:achievable_pfdr} in this manifold setting holds with $\pmin = \pmax = 1/m$ since all target classes are assumed to be equally likely here, and $\dmin = \min_{i \in \{1,\ldots,M\}} d_{i}$ where $$d_i = \min_{\btf_i^{(\ell)},\btf_i^{(k)} \in \sD_{\by_i}, \ell \neq k}\|\btf_i^{(\ell)}-\btf_i^{(k)}\|.$$

\section{Anomalous signal detection}
\label{sec:AnomDet}
The target detection approach discussed above assumes that the target signal of interest resides in a dictionary that is available to the user. However, in some applications (such as military applications and surveillance), one might be interested in detecting objects {\em not} in the dictionary. 
In other words, the target signals of interest are anomalous and are not available to the user. In this section we show how the target detection methods discussed above can be extended to anomaly detection. In particular, we exploit the distance preservation property of the sensing matrix $\bA$ to detect anomalous targets from projection measurements.

\subsection{Problem formulation}
Given observations of the form in \eref{eqn:whitened}, we are interested in
detecting whether $\f^{*} \in \sD$ or $\f^{*}$ is anomalous. Let us
write the anomaly detection problem as the following multiple
hypothesis test:
\begin{subequations}
\begin{align}
\sH_{0i}&:  \|\f^{*}_i-\f\| \leq \tau \text{ for some } \f \in \sD \label{eqn:AnomH0}\\
\sH_{1i}&: \|\f^{*}_i-\f\| > \tau \text{ for all } \f \in \sD \label{eqn:AnomH1}
\end{align}
\end{subequations}
where $\tau \in \left[0,\sqrt{2}\right)$ is a user-defined threshold
that encapsulates our uncertainty about the accuracy with which we
know the dictionary.\footnote{Note that $\tau$ cannot exceed $\sqrt{2}$ because we assume that all targets of interest, including those in $\sD$ and the actual target $\f^*$, are unit-norm.} In particular, $\tau$ controls how different
a signal needs to be from every dictionary element to truly be
considered anomalous. In the absence of any prior knowledge on the
targets of interest, $\tau$ can simply be set to zero. \done{The null hypothesis in this setting models the normal behavior, while the alternative hypothesis models the abnormal or anomalous behavior. This formulation is consistent with the literature \cite{fowler2011anomaly, stein2002adh}.}

Note that the definition of the hypotheses given in \eref{eqn:AnomH0}
and \eref{eqn:AnomH1} matches the definition in \eref{eqn:Hypothesis
  test H1} for the special case where the dictionary contains just one
signal. In this special case, the signal input $\f^{*}$ is
in the dictionary under the null hypothesis in both DSD and ASD problem formulations. {\footnote{The anomaly detection problem discussed here is more accurately described as target detection in the classical detection theory vocabulary. However, in recent works \cite{Hero06geometricentropy,Steinwart05aclassification}, the authors assume that the nominal distribution is obtained from training data and a test sample is declared to be anomalous if it falls outside of the nominal distribution learned form the training data. Our work is in a similar spirit where we learn our dictionary from training data and label any test spectrum that does not correspond to our dictionary as being anomalous.}}


\subsection{Anomaly detection approach}
\label{subsec:anomalyDetApproach}
Our anomaly detection approach and the associated theoretical analysis are based on a ``distance preservation" property of $\bA$, which is stated formally in \eref{eqn:PrePros_distPres}. We propose an anomaly detection method that controls the false discovery rate ($\FDR$) below a desired level $\delta$ for different background and sensor noise statistics. In other words, we control the expected ratio of falsely declared anomalies to the total number of signals declared to be anomalous. 
Note that here we work with the $\FDR$ as opposed to the $\pFDR$, since it is possible for a scene to not contain any anomalies at all. We let $V/R = 0$ for $R = V = 0$ since one does not declare any signal to be anomalous in this case. 
In \cite{benjamini1995cfd}, Benjamini and Hochberg discuss a p-value based procedure, ``BH procedure", that controls the false discovery rate of $M$ independent hypothesis tests below a desired level. Let 
\begin{align}
d_i = \min_{\f \in \sD} \|\by_i-\alpha_i\bA\f\| = \min_{\f \in \sD} \|\alpha_i\bA\left(\f^{*}_i-\f\right) + \bn_i\|
\label{eqn:di}
\end{align}
be the test statistic at the $i\th$ location. The p-value can be defined in terms of our test statistic as follows:
\begin{align}
p_i = \prob\big( \td_i \geq d_i \big| \sH_{0i}\big) \label{eqn:pvalue}
\end{align} 
where ${\td_i} = \min_{\f \in \sD} \|\alpha_i\bA\left(\f^{*}_i-\f\right) +
\bn\|$ and $\bn \sim\sN\left(\zeros,\bI\right)$ is independent of $\bn_i$.
This is the
probability under the null hypothesis, of acquiring a test statistic at least as extreme as the one observed. Let us denote the ordered set
of p-values by $p_{(1)} \leq p_{(2)} \leq \cdots \leq p_{(M)}$ and let
$\sH_{(0i)}$ be the null hypothesis corresponding to $(i)\th$
p-value. The BH procedure says that if we reject all $\sH_{(0i)}$ for
$i = 1,\ldots,t$ where $t$ is the largest $i$ for which $p_{(i)} \leq
i\delta/M$, then the \FDR is controlled at $\delta$. 

To apply this
procedure in our setting, we need to find a tractable expression for
the p-value at every location. This can be accomplished when $\bA$ satisfies the distance-preservation condition stated below. Let $V = \sD \bigcup \{\f^*_i: i \in \{1,\ldots,M\}\}$ be the set of all signals in the dictionary and the ones whose projections are measured. Note that $|V| \leq M+m$.   For a given $\epsilon \in (0,1)$, a projection operator $\bA \in \reals^{K \times N}$, $K \leq N$, is distance-preserving on $V$ if the following holds for all $u,v \in V$:
\begin{align}
(1-\epsilon)\|u - v\| \leq \|\bA(u-v)\| \leq (1+\epsilon)\|u-v\|,  \forall u,v \in V. 
\label{eqn:PrePros_distPres}
\end{align}
{
The existence of such projection operators is guaranteed by the celebrated Johnson and Lindenstrauss (JL) lemma \cite{originalJL}, which says that there exists random constructions of $\bA$ for which \eref{eqn:PrePros_distPres} holds with probability at least $1-2|V|^2 e^{-Kc(\epsilon)}$ provided $K = \Order\left(\log |V|\right) \leq N$, where $c(\epsilon) =\epsilon^2/16-\epsilon^3/48$ \cite{Achlioptas_DatabaseFriendly,baraniuk2008spr}. 
Examples of such constructions are: (a) Gaussian matrices whose entries are drawn from $\sN(0,1/K)$, (b) Bernoulli matrices whose entries are $\pm 1/\sqrt{N}$ with probability $1/2$, (c) random matrices whose entries are $\pm \sqrt{3/N}$ with probability $1/6$ and zero with probability $2/3$ \cite{Achlioptas_DatabaseFriendly,baraniuk2008spr}, and (d) matrices that satisfy the Restricted Isometry Property (RIP) where the signs of the entries in each column are randomized \cite{KW10_2}. }

We now state our main theorem that gives a tight 
upper bound on the p-value at every location when $\{\alpha_i\}$ are unknown and are estimated from the observations. Let $\{\halpha_i\}$ be the estimates of $\{\alpha_i\}$ that satisfy 
\begin{align}
1-\zeta\leq \frac{\alpha_{i}}{\halpha_{i}}\leq 1+\zeta \label{eqn:estAccuracy} 
\end{align}
for $i=1,\ldots,M$ where $\zeta \in [0,1]$ is a measure of the accuracy of the estimation procedure.  
\begin{theorem}
If the $i\th$ hypothesis test is defined according to \eref{eqn:AnomH0} and
\eref{eqn:AnomH1}, the projection matrix $\bA$ satisfies \eref{eqn:PrePros_distPres}
for a given $\epsilon \in (0,1)$, and the estimates $\{\halpha_i\}$ satisfy \eref{eqn:estAccuracy} for some $\zeta \in [0,1]$, then the bound
\begin{align}
p_i &\leq 1-\sF\left(d_i^2; K,(1+\epsilon)^2\halpha_i^2\left(\zeta+\tau\right)^2\right) \label{eqn:pvalueUPalphaEst}
\end{align}
holds for all $i=1,\ldots,M$ 
where $\sF\left(\cdot;
K,\nu\right)$ is the CDF of a noncentral
$\chi^2$ random variable with $K$ degrees of freedom and noncentrality
parameter $\nu$ \cite{wassermanAllofStat}. 
\label{theorem:AnomDet}
\end{theorem}

The proof of this theorem is given in Appendix~\ref{app:AnomDet}. We find the p-value
upper bounds at every location and use the BH procedure to
perform anomaly detection. The performance of this procedure depends
on the values of $K$, $\{\alpha_i\}$, $\tau$ and $\epsilon$. The
parameter $\epsilon$ is a measure of the accuracy with which the
projection matrix $\bA$ preserves the distances between any two
vectors in $\reals^N$. A value of $\epsilon$ close to zero implies
that the distances are preserved fairly accurately. When $\{\alpha_i\}$ are unknown and estimated from the observations, the performance depends on the accuracy of the estimation procedure, which is reflected in our bounds in \eref{eqn:pvalueUPalphaEst} through $\zeta$. 

{One can easily estimate $\{\alpha_i\}$ from $\{\by_i\}$ for some choices of $\bA$. For instance, if the entries of the projection matrix $\bA$ are drawn from $\sN(0,1/K)$,} the $\{\alpha_i\}$ can be estimated using a maximum likelihood estimator (MLE) \done{by exploiting the statistics of the projection matrix and noise. Note that the $\jth$ element of the $\ith$ measured spectrum is $y_{i,j} = \sum_{k = 1}^N \alpha_i f^*_{i,k} a_{j,k} + n_{i,j} \sim \Gaussian{0,\sum_{k = 1}^N \frac{\alpha_i^2}{K} {f^*_{i,k}}^2 + 1}$ for $j \in \setbr{1,\ldots,K}$. Since $\norm{\fstar_i}_2 = 1$ according to our problem formulation, $y_{i,j} \stackrel{\text{i.i.d.}}{\sim} \Gaussian{0,\frac{\alpha_i^2}{K}+1}$. The MLE of $\alpha_i$ given by $\halpha_i = \argmax_{\alpha} \prob(\by_i | \bA, \alpha)$ then reduces to 
\begin{align}
\halpha_i = \sqrt{\paranbr{\|\by_i\|^2-K}}.
\label{eqn:alphahat}
\end{align}
In practice, we use $\halpha_i = \sqrt{\paranbr{\|\by_i\|^2-K}_{+}}$ where the $(a)_{+} = a$ if $a \geq 0$ and $0$ otherwise to ensure that $\|\by_i\|^2-K$ is nonnegative. We can use concentration inequalities to show that with high probability, $\norm{\by_i}_2^2$ is tightly concentrated around its mean $\expect\squarebr{\norm{\by_i}_2^2} = \alpha_i^2+K$. Since $y_{i,j} \stackrel{\text{i.i.d.}}{\sim} \Gaussian{0,\frac{\alpha_i^2}{K}+1}$, $\frac{K}{\alpha^2+K}\norm{\by_i}_2^2 \sim \chi^2_K$. From Lemma 2.2 in \cite{tao2006random}, and Proposition 1 and Remark 1 in \cite{tao_chisq_bounds}, for any $t>0$
\begin{align}
\prob\paranbr{\abs{\norm{\by_i}_2^2 - (\alpha_i^2+K)}\geq t} \leq C \exp(-ct^2)
\end{align}
for some absolute constants $C,c > 0$. This result shows that with high probability, $\norm{\by_i}_2^2-K$ is nonnegative.
}

The experimental results discussed in Sec.~\ref{Sec:ExperimentalResults}
demonstrate the performance of this detector as a function of $K$,
$\{\alpha_i\}$ and $\tau$ when $\{\alpha_i\}$ are known and as a function of $K$, $\tau$ and $\zeta$ when $\{\alpha_i\}$ are estimated.

\section{Experimental Results}
\label{Sec:ExperimentalResults}
{In the experiments that follow, the entries of $\bA$ are drawn from $\sN(0,1/K)$. }
\subsection{Dictionary signal detection}
\label{subsec:DSDExpts}
To test the effectiveness of our approach, we formed a dictionary $\sD$ 
of nine spectra (corresponding to different kinds of trees, grass, water bodies and roads) obtained from a labeled HyMap (Hyperspectral Mapper) remote sensing data set \cite{kruse2000hymap}, and simulated a realistic dataset using the spectra from this dictionary. Each HyMap spectrum is of length $N=106$. We generated projection measurements of these data such that $\bz_i = \alpha_i \bPhi(\f^*_i+\bb_i)+\bw_i$ according
to \eref{eqn:problem setup}, where $\bw_i \sim \sN(0,\sigma^2\bI)$, $\f^*_i \in \sD$ for $i=1,\ldots,8100$, $\bb_{i} \sim \sN\left(\bmu_{\bb},\bSigma_{\bb}\right)$ such that $\bSigma_{\bb}$ satisfies the condition in \eref{eqn:bgdCond}, and $\alpha_{i} = \alpha^*_i\sqrt{K}$ where $\alpha^*_i \sim \sU[21,25]$ and $\sU$ denotes uniform distribution. We let $\sigma^2 = 5$ and model $\{\alpha_i\}$ to be proportional to $\sqrt{K}$ to account for the fact that the total observed signal energy increases as the number of detectors increases.
We transform the $\bz_i$ by a series of
operations to arrive at a model of the form discussed in
\eref{eqn:whitened}, which is $\by_i = \alpha_i \bA\f^*_i + \bn_i$. 
For this dataset, $\pmin = 0.04938$, $\pmax = 0.1481$, and $\dmin = 0.04341$. 

\begin{figure}[hbt]
\centering
\begin{tabular}{cc}
\includegraphics[height=1.75in]{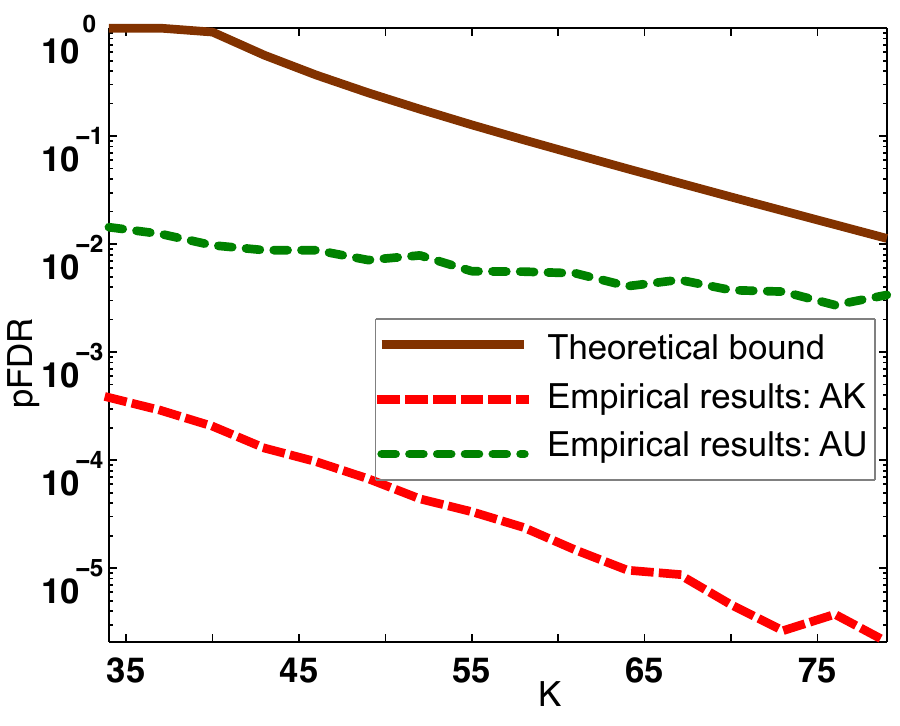}&
\includegraphics[height=1.75in]{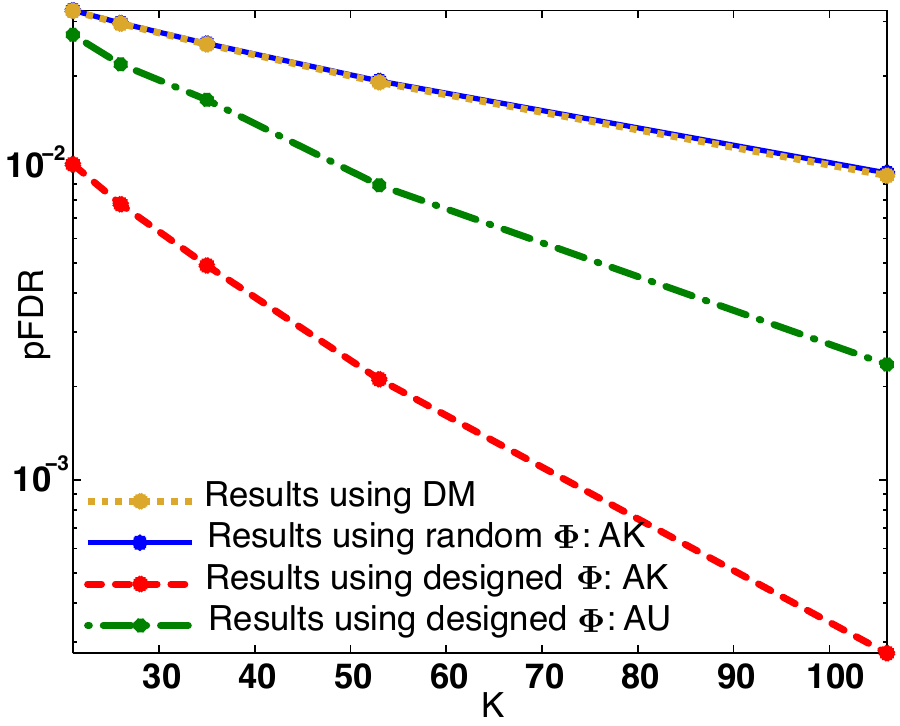}\\
(a) & (b)
\end{tabular}
\small{\caption{Compressive target detection results under the AK ($\setbr{\alpha_i}$ known) and AU ($\setbr{\alpha_i}$ unknown) cases respectively as a function of $K$. (a) Comparison of the worst-case empirical pFDR curves
    with the theoretical bounds when SNR is high. (b) \done{Comparison of the results obtained by the proposed method using projection measurements using $\bPhi$ designed according to \eref{eqn:Phi}, $\bPhi$ chosen at random, and the ones using downsampled measurements (DM) when the SNR is low.}
\label{fig:pFDRCurves}}}
\end{figure}

We evaluate the performance of our detector \eref{eqn:DecisionRule} on the transformed observations, relative to the number of measurements $K$, by comparing the detection results to the ground truth. Our MAP detector returns a label $L^{\text{MAP}}_{i}$ for every observed spectrum which is determined according to 
\begin{eqnarray*}
L^{\text{MAP}}_i = \argmin_{\ell \in \{1,\ldots,m\},\f^{(\ell)} \in \sD} \left(\frac{1}{2}||\by_i - \alpha_i \bA \f^{(\ell)}||^2 - \log p^{(\ell)}\right)
\end{eqnarray*}
where $m$ is the number of signals in $\sD$, and $p^{(\ell)}$ is the a priori probability of target class $\ell$. 
In our experiments we
evaluate the performance of our classifier when (a) $\{\alpha_{i}\}$ are known (AK) and (b) $\{\alpha_i\}$ are unknown (AU) and must be estimated from $\by$, respectively. The empirical $\pFDR^{(j)}$
for each target spectrum $j$ is calculated as follows:
\begin{eqnarray*}
\pFDR^{(j)} = \frac{\sum_{i=1}^{M}\ind{L^{\text{GT}}_i=j}\ind{L^{\text{MAP}}_i\neq j}}{\sum_{i=1}^{M}\ind{L^{\text{MAP}}_i\neq j}}
\end{eqnarray*}
where $\{L^{\text{GT}}_{i}\}$ denote the ground truth labels. The empirical $\pFDR^{(\cdot)}$ is the ratio of the number of missed
targets to the total number of signals that were declared to be
nontargets. The plots in Fig.~\ref{fig:pFDRCurves}(a) show the results obtained using our target detection approach under the AK case (shown by a dark gray dashed line) and the AU case (shown by a light gray dashed line), compared to the theoretical upper bound (shown by a solid line). These results are obtained by averaging the pFDR values obtained over $1000$ different noise, sensing matrix and background realizations. {Note that theoretical results only apply to the AK case since they were derived under the assumption of $\{\alpha_i\}$ being known. The experimental results are shown for both AK and AU cases to provide a comparison between the two scenarios.} In both these cases, the worst-case empirical pFDR curves decay with the increase in the values of $K$. In the AK case, in particular, the worst-case empirical pFDR curve decays at the same rate as the upper bound. {In this experiment, for a fixed $\alpha_{\min}$ and $\dmin$, we chose $K$ to satisfy \eref{eqn:boundonK}. The theory is somewhat conservative, and in practice the method works well even when the values of $K$ are below the bound in \eref{eqn:boundonK}. }

{In the experiment that follows, we let $\alpha^{*}_i \sim \sU[10,20]$, where $\sU$ denotes a uniform random variable, $\alpha_i = \sqrt{K}\alpha^{*}_i$ and evaluate the performance of our detector for different values of $K$ that are not necessarily chosen to satisfy \eref{eqn:boundonK}. In addition, we also compare the performance of our detection method to that of a MAP based target detector operating on downsampled versions of our simulated spectral input image. The reason behind such a comparison is to show what kinds of measurements yield better results given a fixed number of detectors.}

{
For an input spectrum $\bg \in \reals^N$, we let $\btg \in \reals^K$ denote its downsampled approximation. Specifically, the $\jth$ element of $\tg_i$ is $\sum_{\ell = 1}^{r}g_{(j-1)r+\ell}$ where $r = \lceil N/K \rceil$.
Let us consider making observations of the form 
\begin{align}
\by_i = \frac{\btg_i}{c} + \bn_i  \in \reals^K 
\label{eqn:MS_y}
\end{align}
where $\btg_i = \alpha_i\btf^{*}_i + \btb_i$ is the $K$-dimensional downsampled version of $\f^{*}_i + \bb_i$ for $K \leq N$, $\bn_i \sim \sN(\zeros,\sigma^2\bI)$ for $\sigma^2 = 5$ and $c$ is a constant that is chosen to preserve the mean signal-to-noise ratio corresponding to the downsampled and projection measurements. The MAP-based detector operating on the downsampled data returns a label $D^{\text{MAP}}_{i}$ for every observed spectrum which is determined according to 
\begin{align*}
D^{\text{MAP}}_i = \argmin_{\ell \in \{1,\ldots,m\},\f^{(\ell)} \in \sD}\left(\by_i - \alpha_i\btf^{(\ell)}\right)^TG^{-1}&\left(\by_i - \alpha_i\btf^{(\ell)}\right)- \log p^{(\ell)}
\end{align*}
where $G = \btSigma_b + \sigma^2\bI$ and $\btSigma_b$ is the covariance matrix obtained from the downsampled versions of the background training data and $\btf^{(\ell)}$ is the downsampled version of $\f^{(\ell)} \in \sD$. The algorithm declares that target spectrum $\f^{(j)} \in \sD$ is present in the $\ith$ location if $D^{\text{MAP}}_i = j$. 
\done{In order to illustrate the advantages of using a $\bPhi$ designed according to \eref{eqn:Phi}, we compare the performances of the proposed anomaly detector when $\bPhi$ is chosen to be a random Gaussian matrix whose entries are drawn from $\Gaussian{0,1/K}$ and when $\bPhi$ is chosen according to \eref{eqn:Phi}.}
\done{Fig.~\ref{fig:pFDRCurves}(b) shows a comparison of the results obtained using the projection measurements obtained using $\bPhi$ designed according to \eref{eqn:Phi}, $\bPhi$ chosen at random, and the downsampled measurements under the AK case.} These results show that the detection algorithm operating on projection measurements using $\bPhi$ designed using background and sensor noise statistics yield significantly better results than the one operating on the downsampled data, and that the empirical pFDR values in our method decays with $K$. The improvement in performance using projection measurements comes from the distance-preservation property of the projection operator $\bA$. While a Gaussian sensing matrix $\bA$ preserves distances between any pair of vectors from a finite collection of vectors with high probability \cite{Achlioptas_DatabaseFriendly,baraniuk2008spr}, downsampling loses some of the fine differences between similar-looking spectra in the dictionary.} \done{Furthermore, when $\bPhi$ is chosen at random, the resulting whitened transformation matrix is not necessarily distance-preserving. This may worsen the performance as illustrated in Fig.~\ref{fig:pFDRCurves}(b).}

\subsection{Anomaly detection} 
\label{subsec:anomDetExpts}
In this section, we evaluate the performance of our anomaly detection method on (a) a simulated dataset and provide a comparison of the results obtained using the proposed projection measurements and the ones obtained using downsampled measurements, and (b) real AVIRIS (Airborne Visible InfraRed Imaging Spectrometer) dataset. 

\subsubsection{Experiments on simulated data}
We simulate a spectral image $\f^{*}$ composed of $8100$ spectra, where each of them is either drawn from a dictionary $\sD = \{\f^{(1)},\cdots,\f^{(5)}\}$ consisting of five labeled spectra from the HyMap data that correspond to a natural landscape (trees, grass and lakes) or is anomalous. The anomalous spectrum is extracted from unlabeled AVIRIS data, and the minimum distance between the anomalous spectrum $\f^{(\text{a})}$ and any of the spectra in $\sD$ is $\dmin = \min_{\f \in \sD}\|\f-\f^{(\text{a})}\| = 0.5308$. The simulated data has $625$ locations that contain the anomalous spectrum. Our goal is to find the spatial locations that contain the anomalous AVIRIS spectrum given noisy measurements of the form $\bz_i = \bPhi\left(\alpha_i \f^{*}_i + \bb_i\right) + \bw_i$ where $\bb_i \sim (\bmu_b,\bSigma_b)$, $\bPhi$ is designed according to \eref{eqn:Phi}, $\bw_i \sim \sN(\zeros,\sigma^2\bI)$ and $\f^{*}_i \in \sD$ under $\sH_{0i}$. As discussed in Sec.~\ref{sec:AnomDet}, $\f^{*}_i$ is anomalous under $\sH_{1i}$, and our goal is to control the $\FDR$ below a user-specified false discovery level $\delta$. We simulate $\{\alpha_i\} = \sqrt{K}\alpha^{*}_i$ where $\alpha_i^{*} \sim \sU[2,3]$. In this experiment we assume the availability of background training data to estimate the background statistics and the sensor noise variance $\sigma^2$. Given the knowledge of the background statistics, we perform the whitening transformation discussed in Sec.~\ref{sec:preprocessing} and evaluate the detection performance on the preprocessed observations given by \eref{eqn:whitened}. 

\begin{figure}[hbt]
\centering
\begin{tabular}{cc}
\subfigure[Pseudo-ROC plots, GLRT-based method operating on downsampled data using true values of $\setbr{\alpha_i}$]{
\includegraphics[height=1.5in]{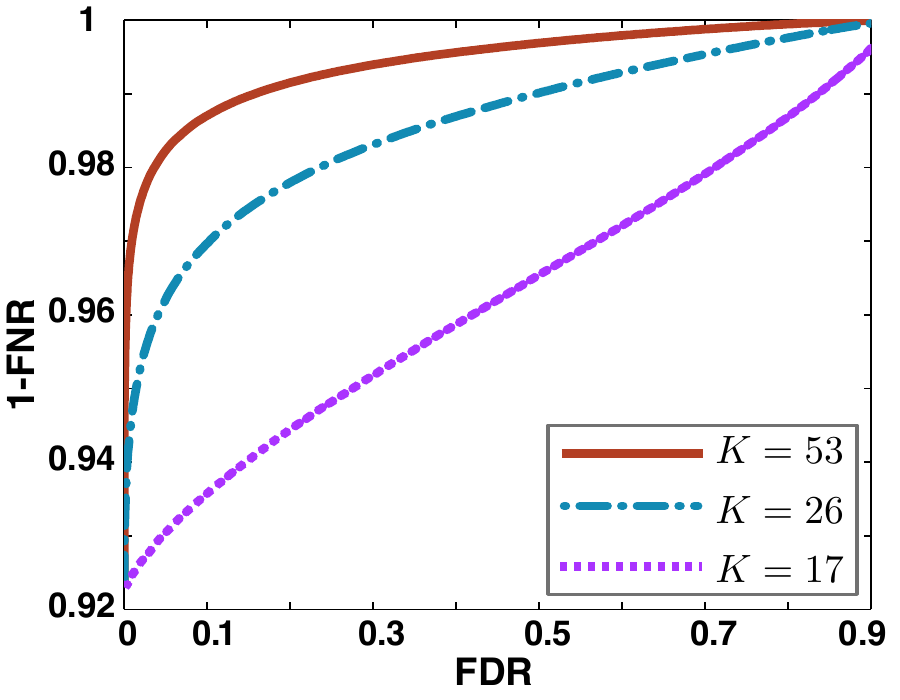}
\label{subfig:FDRFNR_diffK_LRT_AK}
}
\subfigure[Pseudo-ROC plots, Proposed method with $\Phi$ chosen to be a random Gaussian projection matrix using true values of $\setbr{\alpha_i}$]{
\includegraphics[height=1.5in]{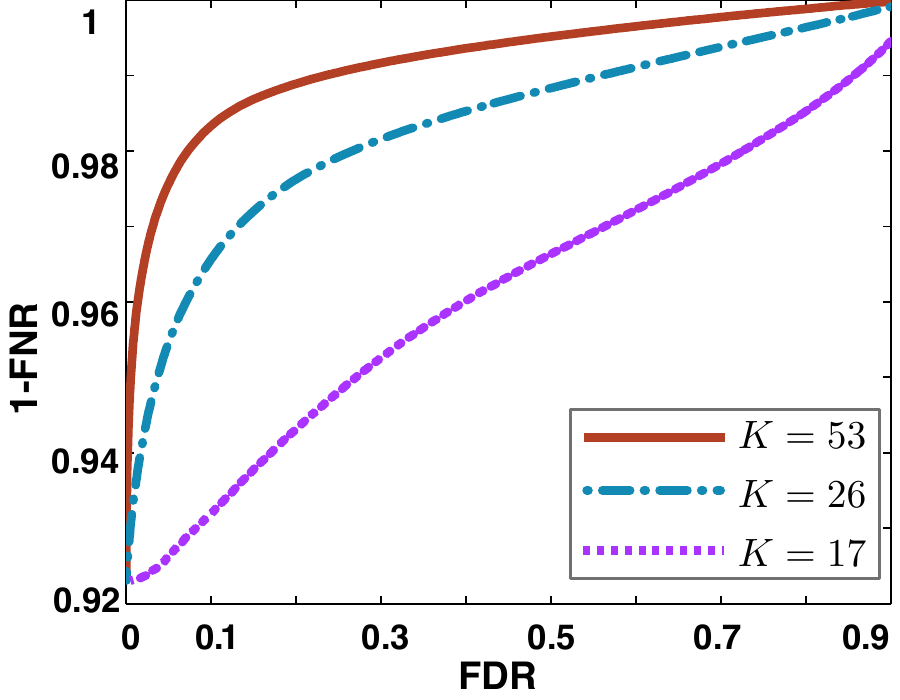}
\label{subfig:FDRFNR_diffK_randProj_AK}
}
\\
\subfigure[Pseudo-ROC plots, Proposed method where $\Phi$ is designed according to \eref{eqn:Phi} using true values of $\setbr{\alpha_i}$]{
\includegraphics[height=1.5in]{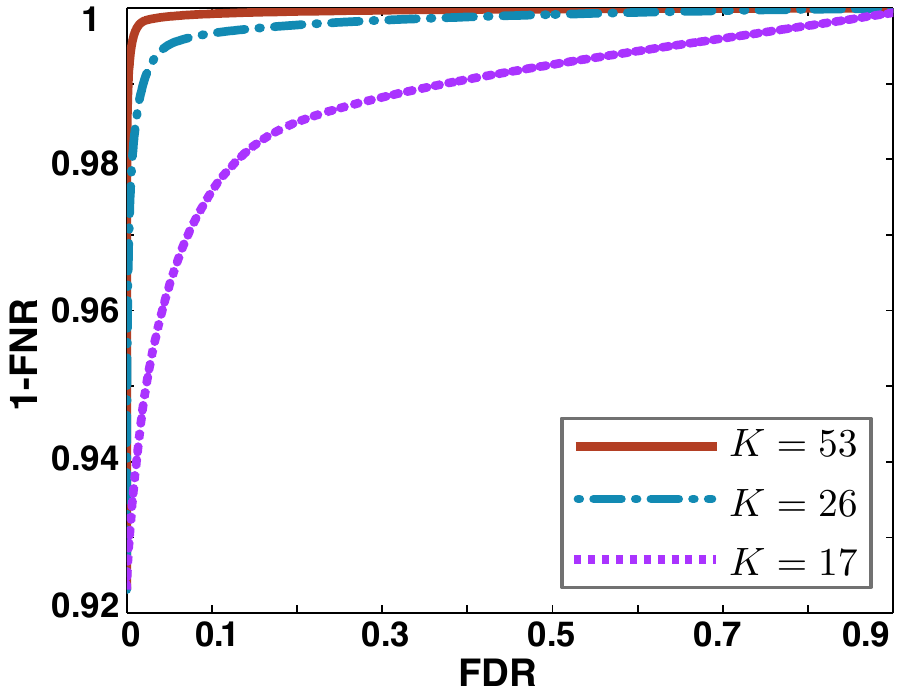}
\label{subfig:FDRFNR_diffK_mine_AK}
}
\subfigure[Pseudo-ROC plots, Proposed method where $\Phi$ is designed according to \eref{eqn:Phi} using ML estimates of $\setbr{\alpha_i}$]{
\includegraphics[height=1.5in]{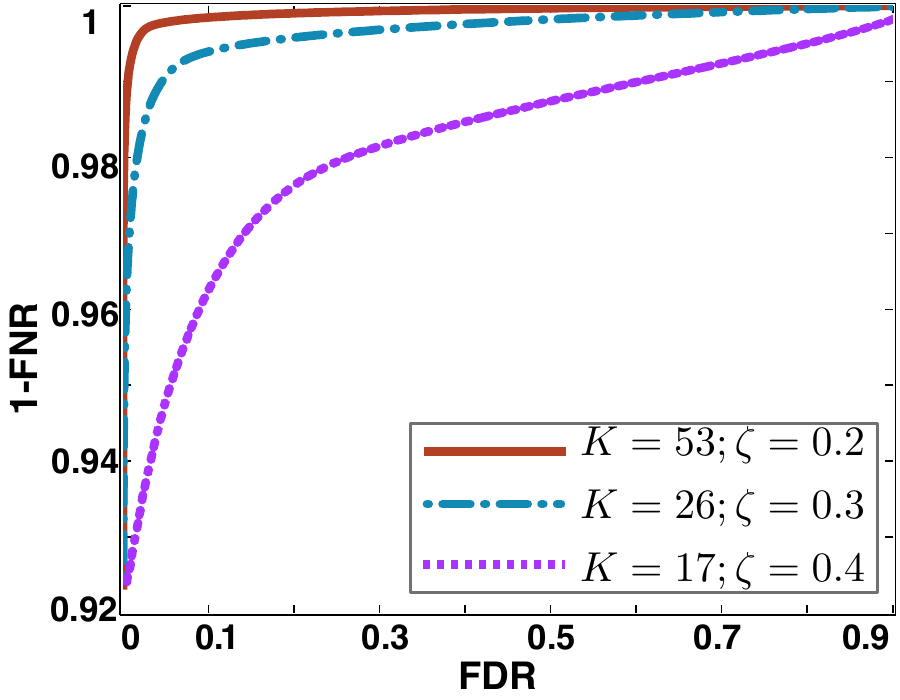}
\label{subfig:FDRFNR_diffK_mine_AU}
}
\\
\subfigure[ROC plots, GLRT-based method operating on downsampled data using true values of $\setbr{\alpha_i}$]{
\includegraphics[height=1.5in]{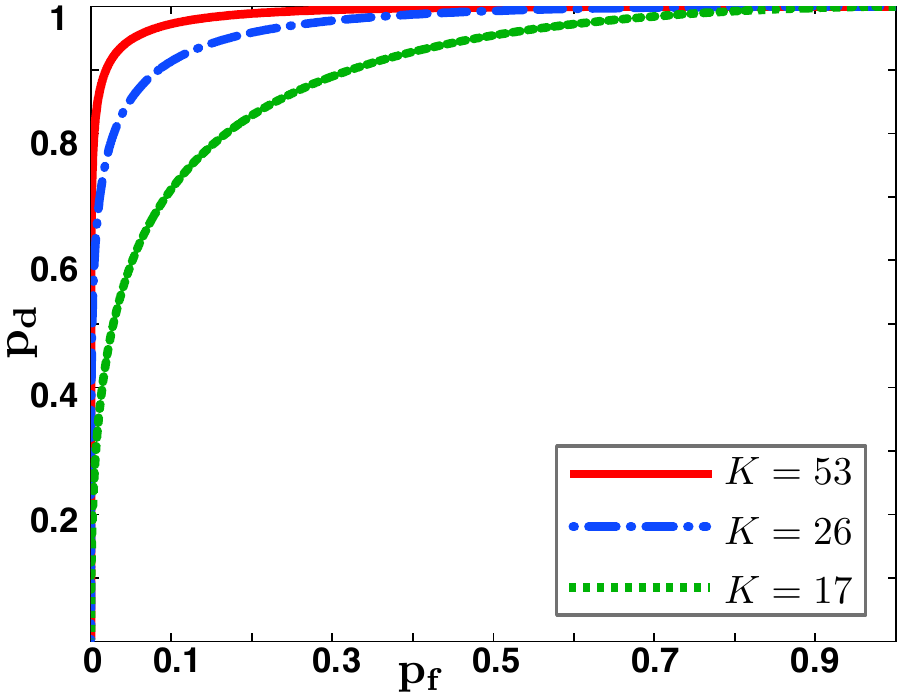}
\label{subfig:pdpf_diffK_LRT_AK}
}
\subfigure[ROC plots, Proposed method with $\Phi$ chosen to be a random Gaussian projection matrix using true values of $\setbr{\alpha_i}$]{
\includegraphics[height=1.5in]{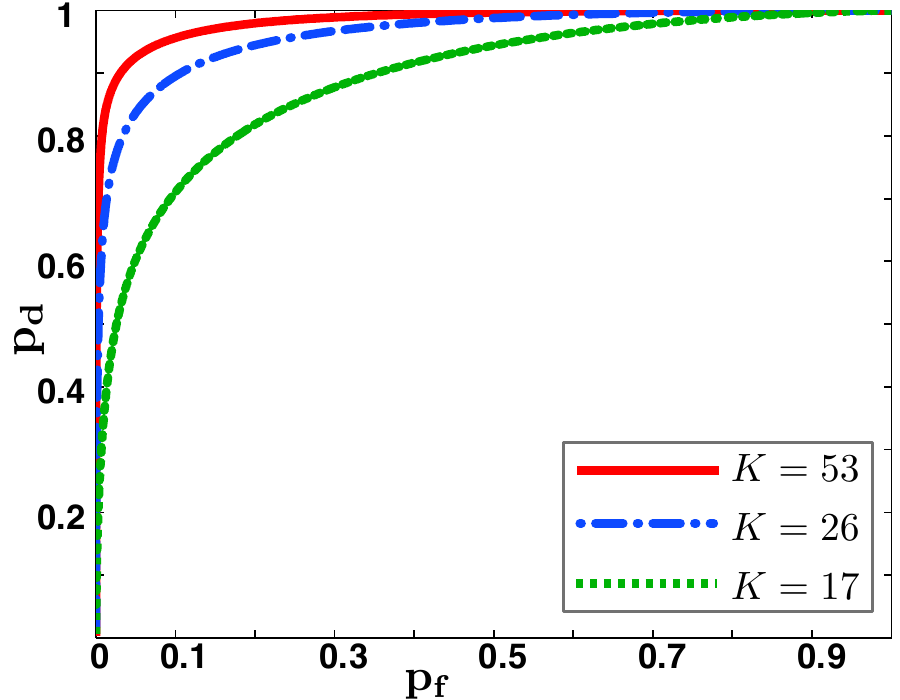}
\label{subfig:pdpf_diffK_randProj_AK}
}\\
\subfigure[ROC plots, Proposed method where $\Phi$ is designed according to \eref{eqn:Phi} using true values of $\setbr{\alpha_i}$]{
\includegraphics[height=1.5in]{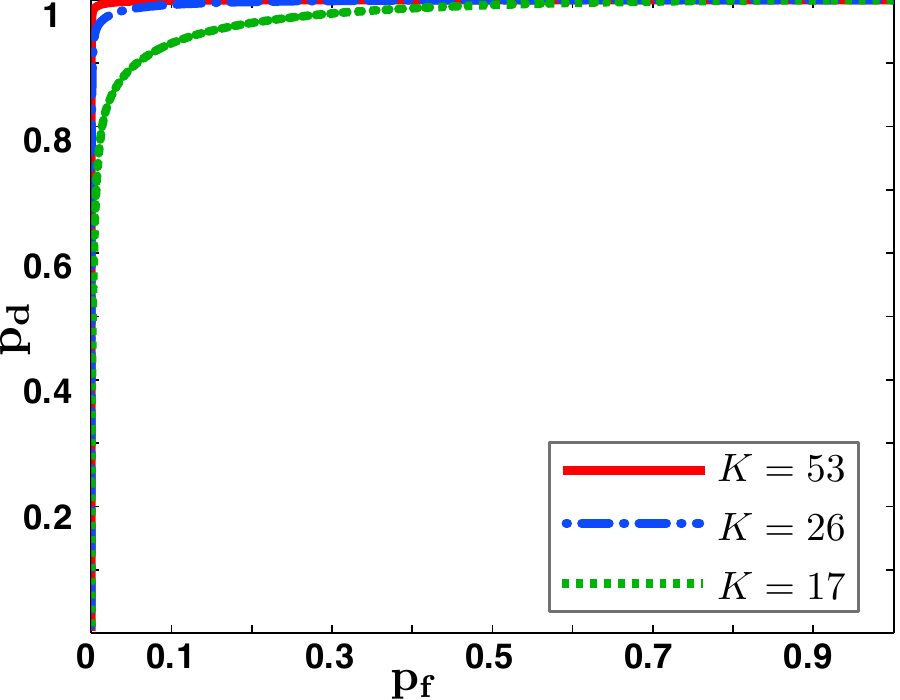}
\label{subfig:pdpf_diffK_mine_AK}
}
\subfigure[ROC plots, Proposed method where $\Phi$ is designed according to \eref{eqn:Phi} using ML estimates of $\setbr{\alpha_i}$]{
\includegraphics[height=1.5in]{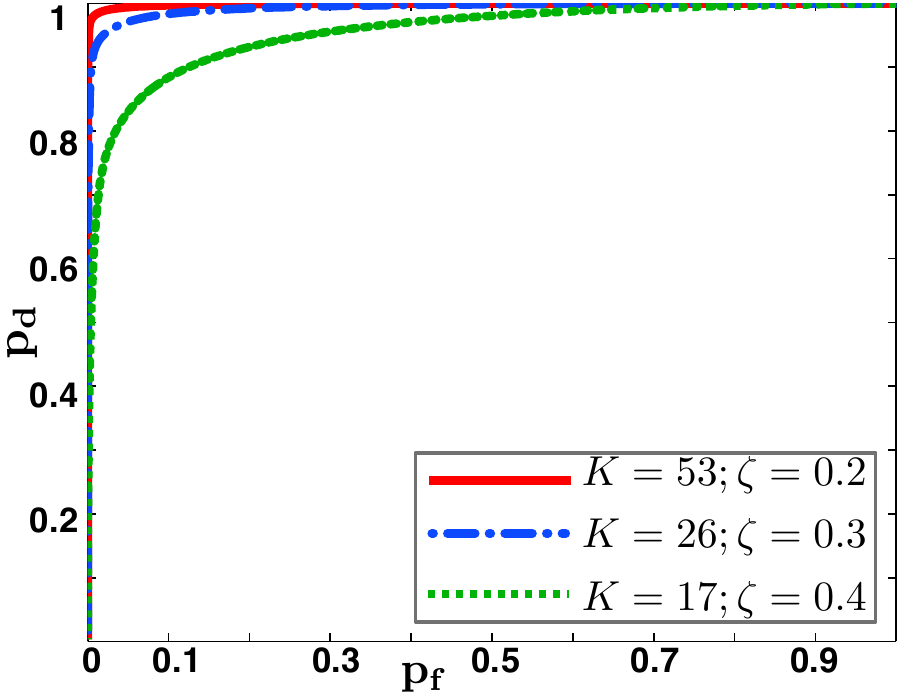}
\label{subfig:pdpf_diffK_mine_AU}}
\end{tabular}
\small{\caption{Comparison of the performances of the proposed anomaly detector using a random $\Phi$, the proposed anomaly detector using the designed $\Phi$ in \eref{eqn:Phi} and the GLRT-based method operating on downsampled data for different values of $K$ when $\alpha_i^* \in \sU[2,3]$ and $\alpha_i = \alpha_i^*\sqrt{K}$.
\label{fig:FDRFNR_diffK}}}
\end{figure}
For a fixed $\tau = 0.1$ and $\epsilon = 0.1$, we evaluate the performance of the detector as the number of measurements $K$ increases under the AK and AU cases respectively, by comparing the \emph{pseudo-ROC} (receiver operating characteristic) curves obtained by plotting the empirical false discovery rate against $1-\FNR$, where $\FNR$ is the false nondiscovery rate. Note that $1-\FNR$ is the expected ratio of the number of null hypotheses that are correctly rejected to the number of declared null hypotheses. The empirical $\FDR$ and $\FNR$ are computed according to
\begin{align*}
\FDR =   \frac{\sum_{i=1}^M\identity_{\left\{L^{\text{GT}}_i=0\right\}}\identity_{\{p_i \leq p_t\}}}{\sum_{i=1}^M\identity_{\{p_i \leq p_t\}}} \text{ and }
\FNR =  \frac{\sum_{i=1}^M \identity_{\left\{L^{\text{GT}}_i=1\right\}}\identity_{\{p_i > p_t\}}}{\sum_{i=1}^M \identity_{\{p_i > p_t\}}} 
\end{align*}
where $p_t$ is the p-value threshold such that the BH procedure rejects all null hypotheses for which $p_i \leq p_t$, and the ground truth label $L^{\text{GT}}_i = 0$ if the $\ith$ spectrum is not anomalous, and $1$ otherwise. In this experiment, we consider three different values of $K$ approximately given by $K \in \{N/6, N/3, N/2\}$ where $N=106$, and evaluate the performance of our detector for each $K$. Furthermore, in our experiments with simulated data, we declare a spectrum to be anomalous if $d_i \geq \eta$ where $\eta$ is a user-specified threshold and $d_i$ is defined in \eref{eqn:di}. We use the p-value upper bound in \eref{eqn:pvalueUPalphaEst} in our experiments with real data where the ground truth is unknown. 

{We compare the performance of our method to a generalized likelihood ratio test (GLRT)-based procedure operating on downsampled data, where we collect measurements of the form in \eref{eqn:MS_y} and $\f^{*}_i \in \sD$ under $\sH_{0i}$. Observe that $\by_i|\sH_{0i} \sim \sum_{\f \in \sD}\prob\left(\f^{*}_i = \f\right)\sN(\alpha_i\btf,\btSigma_b + \bI)$, where $\btf$ refers to the downsampled version of $\f \in \sD$. In this experiment we assume that each spectrum in $\sD$ is equally likely under $\sH_{0i}$ for $i=1,\ldots,M$. The GLRT-based approach declares the $\ith$ spectrum to be anomalous if 
\begin{align*}
-\log \prob\left(\by_i | \sH_{0i}\right) \overset{\sH_{1i}}{\underset{\sH_{0i}}{\gtrless}} \eta 
\end{align*}
for $i=1,\ldots,M$, where $\eta$ is a user-specified threshold \cite{stein2002adh}. While our anomaly detection method is designed to control the $\FDR$ below a user-specified threshold, the GLRT-based method is designed to increase the probability of detection while keeping the probability of false alarm as low as possible. To facilitate a fair evaluation of these methods, we compare the pseudo-ROC curves ($\FDR$ versus $1-\FNR$) and the actual ROC curves (probability of false alarm $p_f$ versus probability of detection $p_d$) corresponding to these methods obtained by averaging the empirical $\FDR$, $\FNR$, $p_d$ and $p_f$ over $1000$ different noise and sensing matrix realizations for different values of $K$. 
\done{We also compare the performance of the proposed method when $\bPhi$ is chosen according to \eref{eqn:Phi} and when it is chosen at random, as discussed in the previous section.}
Figs.~\ref{fig:FDRFNR_diffK}\subref{subfig:FDRFNR_diffK_LRT_AK} and \ref{fig:FDRFNR_diffK}\subref{subfig:pdpf_diffK_LRT_AK} show the pseudo-ROC plots and the conventional ROC plots obtained using the GLRT-based method operating on downsampled data when $\{\alpha_i\}$ are known. 
\done{Figs.~\ref{fig:FDRFNR_diffK}\subref{subfig:FDRFNR_diffK_randProj_AK} and \ref{fig:FDRFNR_diffK}\subref{subfig:pdpf_diffK_randProj_AK} show the results obtained by using a random Gaussian $\bPhi$ instead of the $\bPhi$ in \eref{eqn:Phi}. }
Figs.~\ref{fig:FDRFNR_diffK}\subref{subfig:FDRFNR_diffK_mine_AK} and \ref{fig:FDRFNR_diffK}\subref{subfig:pdpf_diffK_mine_AK} show the pseudo-ROC plots and the conventional ROC plots obtained using our method when $\{\alpha_i\}$ are known. These plots show that performing anomaly detection from \done{our \emph{designed}} projection measurements yields better results than performing anomaly detection on downsampled measurements and \done{on measurements obtained using a random Gaussian $\bPhi$}. 
This is largely due to the fact that carefully chosen projection measurements preserve distances (up to a constant factor) among pairs of vectors in a finite collection, where as the downsampled measurements fail to preserve distances among vectors that are very similar to each other. 
\done{Similarly, a random projection matrix $\bPhi$ is not necessarily distance-preserving post-whitening transformation, which leads to poor performance as illustrated in Figs.~\ref{fig:FDRFNR_diffK}\subref{subfig:FDRFNR_diffK_randProj_AK} and \ref{fig:FDRFNR_diffK}\subref{subfig:pdpf_diffK_randProj_AK}.} 
Figs.~\ref{fig:FDRFNR_diffK}\subref{subfig:FDRFNR_diffK_mine_AU} and \ref{fig:FDRFNR_diffK}\subref{subfig:pdpf_diffK_mine_AU} show the pseudo-ROC plots and the conventional ROC plots obtained using our method when $\{\alpha_i\}$ are unknown, and are estimated from the measurements. Note that the value of $\zeta$ decreases as $K$ increases since the estimation accuracy of $\{\alpha_i\}$ increases with increase in $K$. These plots show that the performance improves as we collect more observations, and that, as expected, the performance under the AK case is better than the performance under the AU case. 
{
\subsubsection{Experiments on real AVIRIS data}
To test the performance of our anomaly detector on a real dataset, we consider the unlabeled AVIRIS Jasper Ridge dataset $\bg \in \reals^{614 \times 512 \times 197}$, which is publicly available from the NASA AVIRIS website, \url{http://aviris.jpl.nasa.gov/html/aviris.freedata.html}. We split this data spatially to form equisized training and validation datasets, $\bg^{t}$ and $\bg^v$ respectively, each of which is of size $128 \times 128 \times 197$. Figs.~\ref{fig:AVIRISresults}(a) and \ref{fig:AVIRISresults}(b) show images of the AVIRIS training and validation data summed through the spectral coordinates. The training data are comprised of a rocky terrain with a small patch of trees. The validation data seems to be made of a similar rocky terrain, but also contain an anomalous lake-like structure. The goal is to evaluate the performance of the detector in detecting the anomalous region in the validation data for different values of $K$. We cluster the spectral targets in the normalized training data to eight different clusters using the K-means clustering algorithm and form a dictionary $\sD$ comprising of the cluster centroids. Given the dictionary and the validation data, we find the ground truth by labeling the $\ith$ validation spectrum as anomalous if $\min_{\f \in \sD}\left\|\f-\frac{\bg_i^v}{\|\bg_i^v\|}\right\| > \tau$. Since the statistics of the possible background contamination in the data could not be learned in this experiment because of the lack of labeled training data, the dictionary might be background contaminated as well. The parameter $\tau$ encapsulates this uncertainty in our knowledge of the dictionary. In this experiment, we set $\tau = 0.2$. }

{We generate measurements of the form $\by_i = \sqrt{K}\bg^v_i + \bn_i$ for $i = 1,\ldots,128\times 128$, where $\bn_i \sim \sN(\zeros,\bI)$. The $\sqrt{K}$ factor indicates that the observed signal strength increases with $K$. For a fixed $\FDR$ control value of $0.01$, Figs.~\ref{fig:AVIRISresults}(c) and \ref{fig:AVIRISresults}(d) show the results obtained for $K \approx N/5$ and $K \approx N/2$ respectively. Fig.~\ref{fig:AVIRISresults}(e) shows how the probability of error decays as a function of the number of measurements $K$. The results presented here are obtained by averaging over $1000$ different noise and sensing matrix realizations. From these results, we can see that the number of detected anomalies increases with $K$ and the number of misclassifications decrease with $K$. }

\begin{figure}[hbt]
\centering
\subfigure[Training data]{
\label{subfig:training}
\includegraphics[height=1.25in]{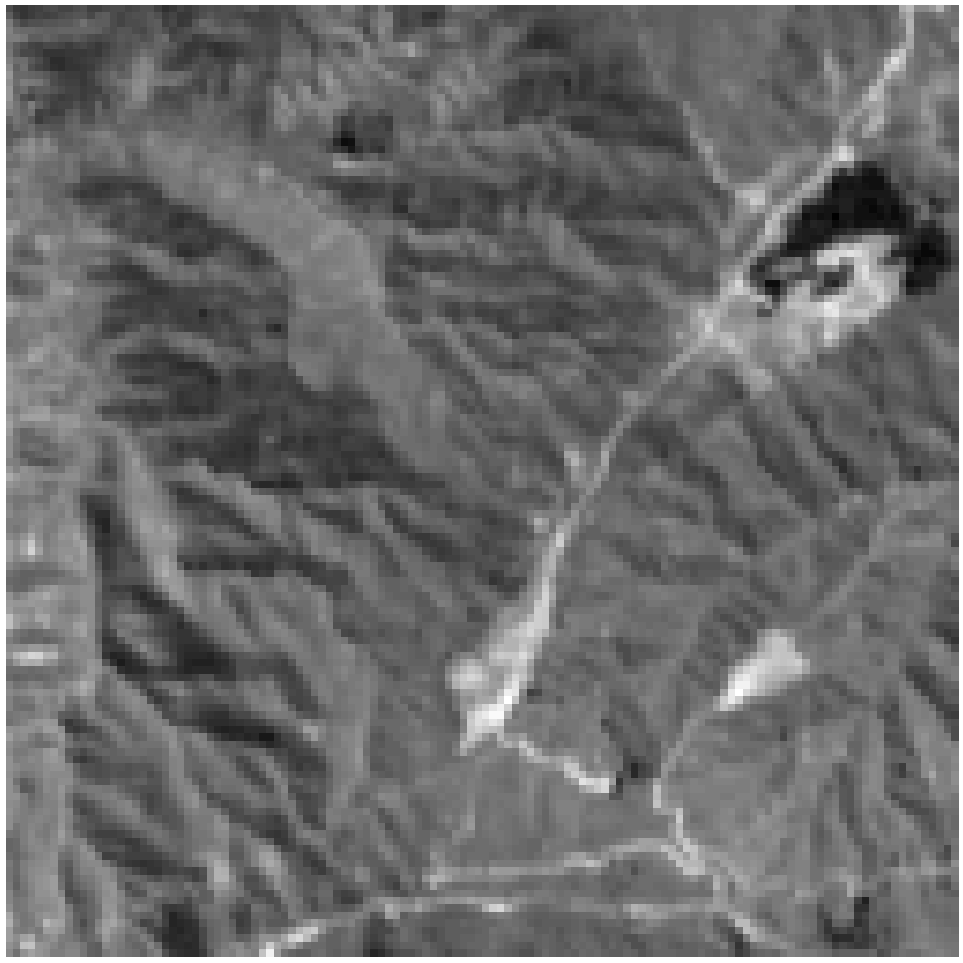}
}
\subfigure[Validation data]{
\includegraphics[height=1.25in]{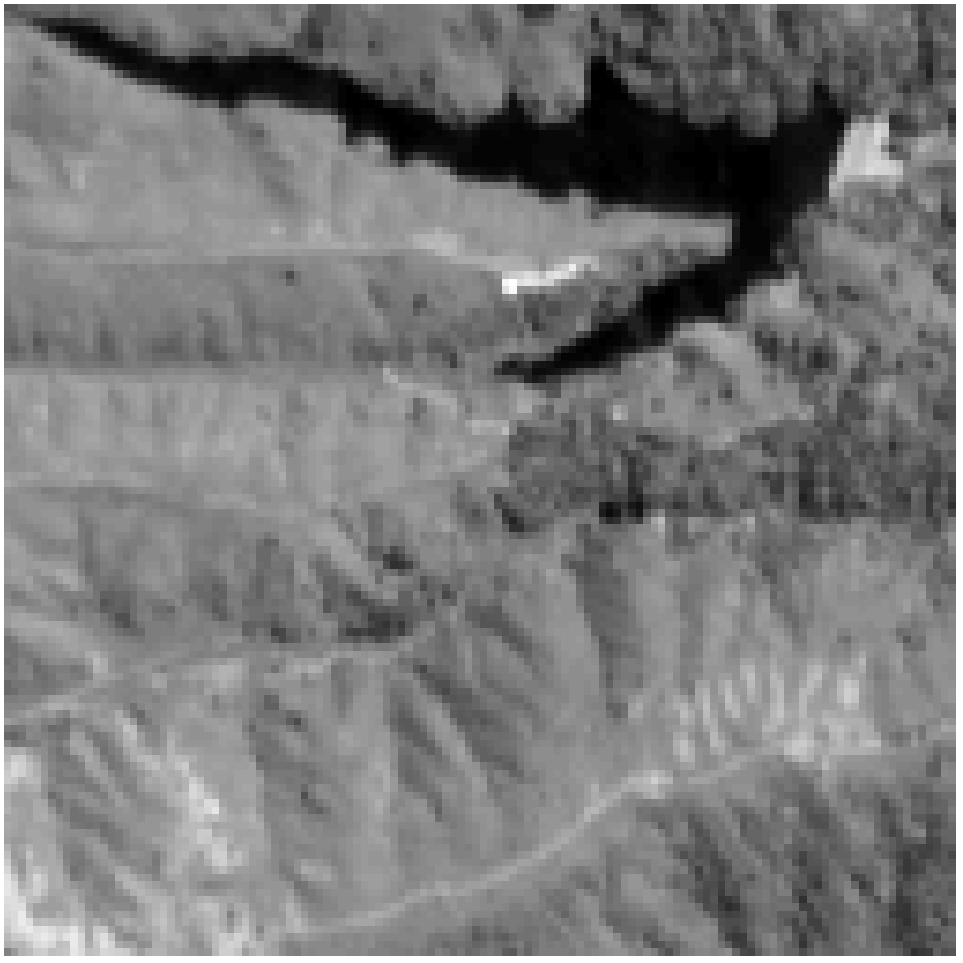}
}
\subfigure[Anomalies detected (shown by white dots) for $K \approx N/5 = 39$.]{
\includegraphics[height=1.25in]{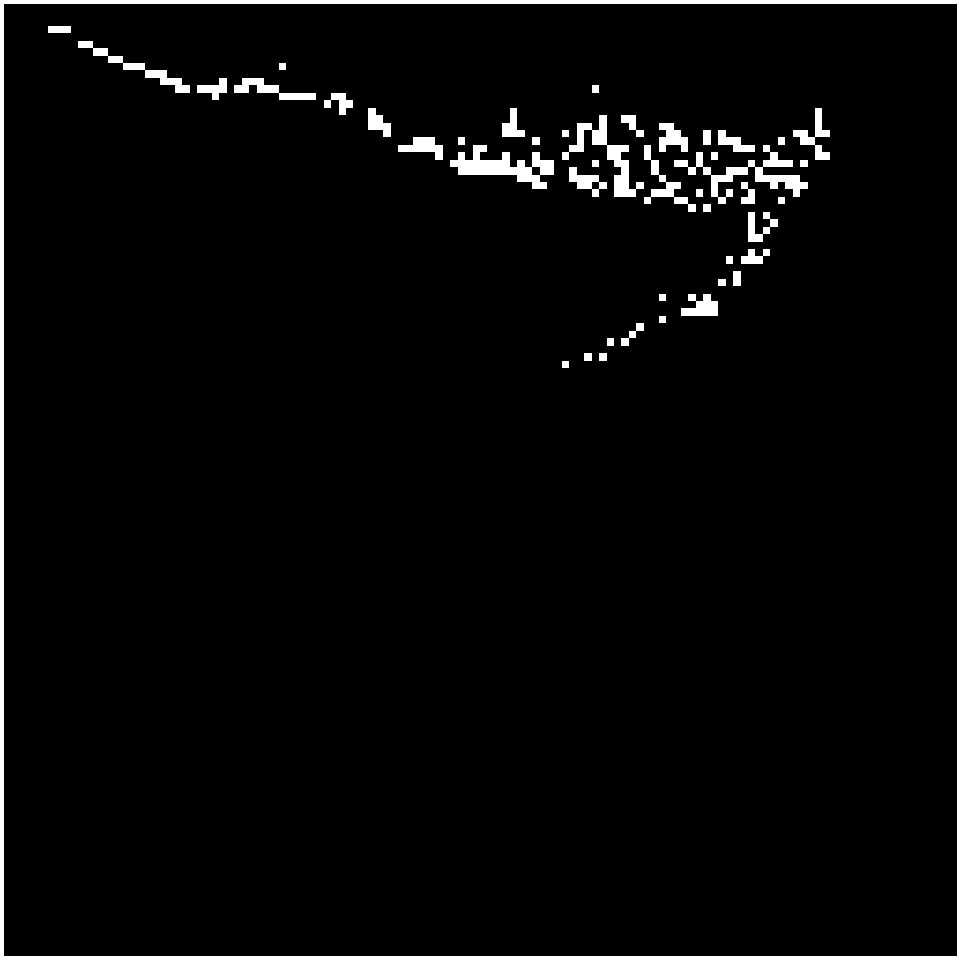}
}
\subfigure[Anomalies detected (shown by white dots) for $K \approx N/2 = 99$.]{
\includegraphics[height=1.25in]{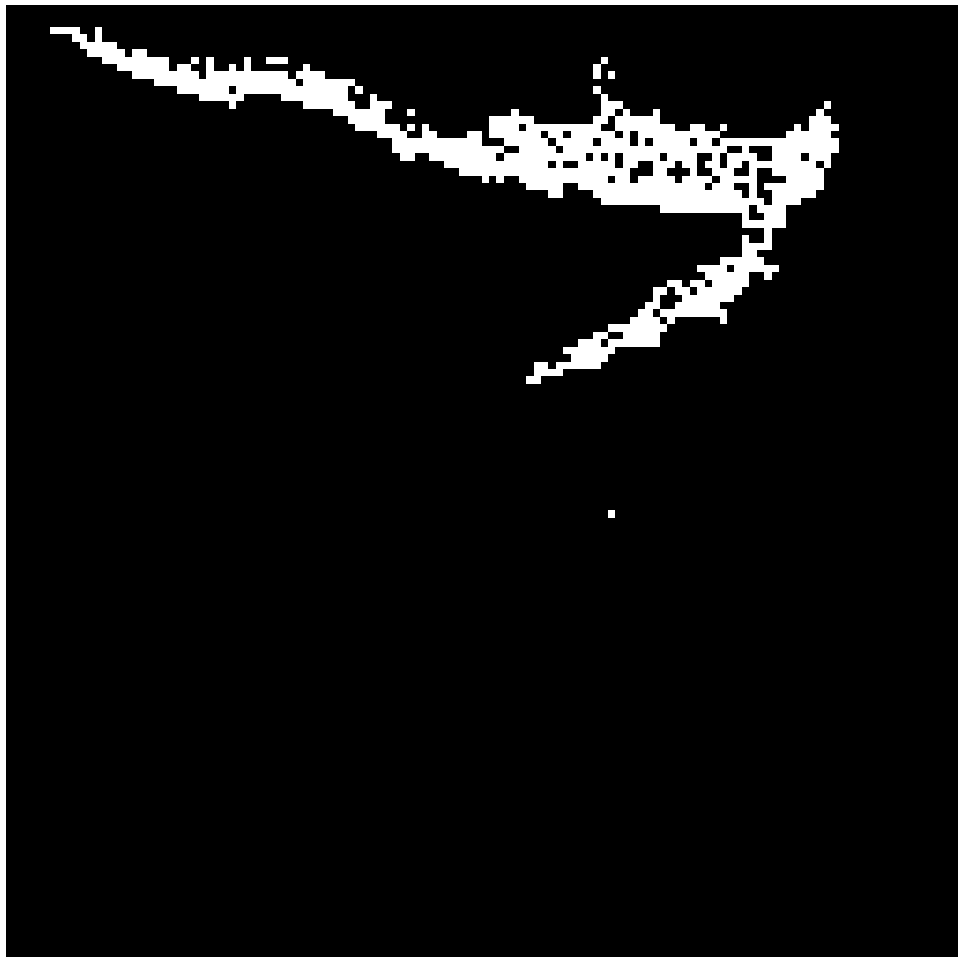}
}
\subfigure[Plot of the probability of error $p_e$ for different values of $K$.]{{\includegraphics[height=1.25in]{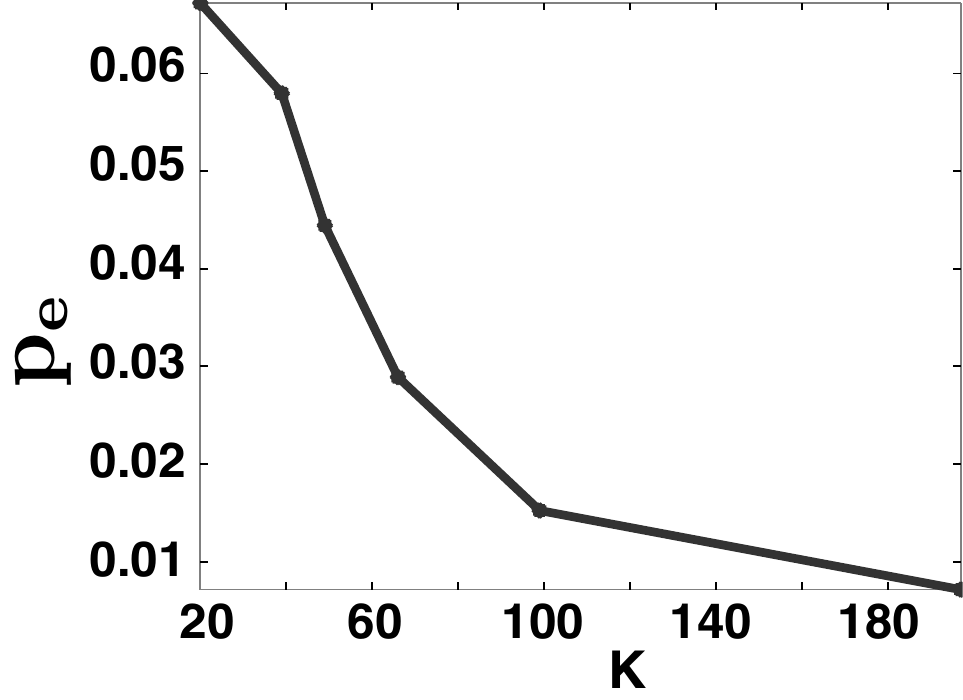}}
}
\caption[Anomaly detection results corresponding to real AVIRIS data.]{\small{Anomaly detection results corresponding to real AVIRIS data for a fixed $\FDR$ control of $0.01$. }}
\label{fig:AVIRISresults}
\end{figure}


\section{Conclusion}
\label{sec:conclusion}
This work presents computationally efficient approaches for
detecting known targets and anomalies of different strengths from
projection measurements without performing a complete
reconstruction of the underlying signals, and offers theoretical bounds
on the worst-case target detector performance.
This paper treats each signal as independent
of its spatial or temporal neighbors. This assumption is reasonable in many
contexts, especially when the spatial or temporal resolution is low relative to
the spatial homogeneity of the environment or the pace with which a scene changes. However, emerging
technologies in computational optical systems continue to improve the
resolution of spectral imagers. In our future work we will build upon
the methods that we have discussed here to exploit the spatial or temporal
correlations in the data.

\begin{appendix}
\section{Proof of Theorem~\ref{thm:PhiConstruction}}
\label{app:PhiConstruction}
Using linear algebra and matrix theory, it is possible to show that if $\bB = \bI - \bA \bSigma_b \bA^T$ is positive definite, then 
\begin{align}
\bPhi = \sigma \bB^{-1/2} \bA \label{eqn:Phi}
\end{align}
satisfies \eref{eqn:CPhiPhi}.\footnote{The authors would like to thank Prof.~Roummel Marcia for fruitful discussions related to this point.} In particular, we can substitute \eref{eqn:Phi} in \eref{eqn:CPhiPhi} to verify that the proposed construction of $\bPhi$ satisfies \eref{eqn:CPhiPhi}. Observe that $C_{\bPhi} = \left(\bPhi \bSigma_b \bPhi^T + \sigma^2\bI \right)^{-1/2}$ can be written in terms of \eref{eqn:Phi} as follows:
\begin{align}
C_{\bPhi}&= \left(\left[\sigma \bB^{-\frac{1}{2}} \bA\right] \bSigma_b \left[\sigma \bB^{-\frac{1}{2}} \bA\right]^T + \sigma^2\bI \right)^{-\frac{1}{2}} \nonumber \\
&= \left(\sigma^{2} \bB^{-1/2} \left(\bA \bSigma_b \bA^T\right)\left(\bB^{-\frac{1}{2}}\right)^{T}  + \sigma^2\bI \right)^{-\frac{1}{2}} \nonumber \\
&=  \left(\sigma^{2} \bB^{-\frac{1}{2}} \left(\bI-\bB\right)\left(\bB^{-\frac{1}{2}}\right)^{T}  + \sigma^2\bI \right)^{-\frac{1}{2}} 
= \left(\sigma^{2} \bB^{-1}\right)^{-\frac{1}{2}} = \sigma^{-1} \bB^{\frac{1}{2}} \label{eqn:last} 
\end{align}
where the third-to-last equation follows from the definition of $\bB$ and \eref{eqn:last} follows from the fact that $\bB$ is symmetric and positive definite. If $\bB$ is positive definite, then $\bB^{-1}$ is positive definite as well and can be decomposed as $\bB^{-1} = \left(\bB^{-1/2}\right)^{T}\bB^{-1/2}$, where the matrix square root $\bB^{-1/2}$ is symmetric and positive definite. By substituting \eref{eqn:last} and \eref{eqn:Phi} in \eref{eqn:CPhiPhi}, we have $C_{\bPhi}\bPhi = \sigma^{-1} \bB^{1/2} \sigma \bB^{-1/2} \bA = \bA$. A sufficient condition for $\bB$ to be positive definite can be derived as follows.

To ensure positive definiteness of $\bB$, we must have
\begin{align}
\bx^{T}\bB\bx &= \bx^{T}\bx - \bx^T\left(\bA \bSigma_{\bb} \bA^T\right)\bx > 0 \label{eqn:lemma1eqn}
\end{align}
 for any nonzero $\bx \in \reals^K$. Note that since $\bSigma_{\bb}$ is positive semidefinite, $\bx^T\left(\bA \bSigma_{\bb} \bA^T\right)\bx \geq 0$. However, the right hand side of \eref{eqn:lemma1eqn} is $>0$ only if the spectral norm of $\bA \bSigma_{\bb} \bA^T$ is $<1$, since $\bx^T\left(\bA \bSigma_{\bb} \bA^T\right)\bx \leq \|\bx\|^2 \cdot \|\bA \bSigma_{\bb} \bA^T\|$. The norm of $\bA \bSigma_{\bb} \bA^T$ is in turn bounded above by
 \begin{align*}
 \|\bA \bSigma_{\bb} \bA^T\| &\leq \|\bA\| \|\bSigma_{\bb}\| \|\bA^{T}\| = \|\bA\|^2 \|\bSigma_{\bb}\| = \|\bA\|^2 \lambda_{\max} 
 \end{align*}
since $\|\bA\| = \|\bA^{T}\|$ and $\|\bSigma_{\bb}\| = \lambda_{\max}$, where $\lambda_{\max}$ is the largest eigenvalue of $\bSigma_{\bb}$. To ensure $\|\bA \bSigma_{\bb} \bA^T\| < 1$, $\|\bA\|^2 \lambda_{\max}$ has to be $<1$, which leads to the result of Theorem~\ref{thm:PhiConstruction}. 
\section{Proof of Theorem~\ref{thm:pFDR}} 
\label{app:pFDR}
The proof of Theorem~\ref{thm:pFDR} adapts the proof techniques from \cite{storey2003positive} to nonidentical independent hypothesis tests. 
We begin by expanding the pFDR definition in \eref{eqn:pFDRdefn} as follows:
\begin{align*}
\pFDRj\bGammabr = \sum_{k=1}^{M}  \expect &\left[\left.\frac{V\bGammabr}{R\bGammabr} \right| R\bGammabr =k\right] \times \prob\left(\left.R\bGammabr = k \right| R\bGammabr > 0 \right). 
\end{align*}
Observe that $R\bGammabr = k$ implies that there exists some subset $S_k = \{u_1,\ldots,u_k\} \subseteq \{1,\ldots,M\}$ of size $k$ such that{$\by_{u_\l} \in \Gamma_{u_\l}^{(j)}$ for $\l = 1,\ldots,k$} and $\by_{i} \not\in \Gamma_i^{(j)}$ for all $i \not\in S_k$. To simplify the notation, let $\Lambda_{S_k} = \prod_{u \in S_k} \Gamma_u^{j} \times \prod_{\ell \notin S_k} \tGamma^{(j)}_{\ell}$, where $\tGamma^{(j)}_{\ell}$ is the complement of $\Gamma^{(j)}_{\ell}$, denote the significance region that corresponds to set $S_k$, and $\bT = (\by_1, \ldots,\by_{M})$ be a set of test statistics corresponding to each hypothesis test. Considering all such subsets we have
\begin{align}
\pFDRj\bGammabr = \sum_{k=1}^{M} \sum_{S_k} \expect& \left[\left.\frac{V\bGammabr}{k} \right| \bT \in \Lambda_{S_k}\right] \times\prob \left(\left.\bT \in \Lambda_{S_k} \right| R\bGammabr > 0 \right). 
\label{eqn: pFDRSubsetEqn}
\end{align}
By plugging in the definition of $V\setGamma$ from \eref{eqn:Vdefn}, we have
\begin{align}
\expect &\left[\left.V\bGammabr \right| \bT \in \Lambda_{S_k}\right]  = \expect \left[\left.\sum_{i=1}^{M} \ind{\by_{i} \in \Gamma_i^{(j)}} \ind{\sH_{i}^{(j)}=0} \right| \bT \in \Lambda_{S_k} \right] \nonumber \\
&\equiv \sum_{\l=1}^{k} \expect \left[ \left.\ind{\sH_{u_\l}^{(j)}=0}\right|\by_{u_\l} \right] = 
\sum_{\l=1}^{k}\prob\left(\left.\sH_{u_\l}^{(j)}=0 \right| \by_{u_\l} \in \Gamma_{u_\l}^{(j)} \right) \label{eqn:expectation_expn}
\end{align}
for all $u_\l \in S_k$ since the tests are independent of each other given $\bA$. 
The posterior probability $\prob\left(\left.\sH_{i}^{(j)}=0 \right| \by_{i} \in \Gamma_{i}^{(j)} \right)$ for the $\ith$ hypothesis test can be expanded using Bayes' rule as 
\begin{align}
\prob\left(\sH_{0i}^{(j)}\left|\by_i \in \Gamma_i^{(j)}\right.\right) &= \frac{\prob\left(\by_i \in \Gamma_i^{(j)}\left|\sH_{0i}\right.\right)\prob\left(\sH_{0i}^{(j)}\right)}{\prob\left(\by_i^{(j)} \in \Gamma_i^{(j)}\right)}\nonumber \\
&\equiv \frac{\prob\left(\left.\bhf_i \neq \f^{(j)}\right|\f^{*}_i = \f^{(j)}\right)\prob\left(\f^{*}_i = \f^{(j)}\right)}{\prob\left(\bhf_i \neq \f^{(j)}\right)},
\label{eqn:pFDR_oneTest}
\end{align}
{where $\bhf_i = \argmax_{\f^{(\ell)} \in \sD}\prob\left(\left.\f_i^{*}=\f^{(\ell)}\right|\by_i,\alpha_{i},\bA\right)$. To upper bound the numerator of \eref{eqn:pFDR_oneTest}, consider the probability of misclassification given by $\pe_{i}= \prob\left(\bhf_i \neq \f^{*}_i\right)$ where $\f^{*}_i = \f^{(j)} \in \sD$, which can be expanded as follows:
\begin{align}
&\pe_i = \prob\left(\bhf_i \neq \f^{*}_i\right) = \sum_{\ell=1}^{m}\prob\left(\left.\bhf_i \neq \f^{*}_i\right|\f^{*}_i = \f^{(\ell)}\right)\prob\left(\f^{*}_i=\f^{(\ell)}\right)\nonumber\\
&\equiv \sum_{\ell=1}^{m}\prob\left(\left.\bhf_i \neq \f^{(\ell)}\right|\f^{*}_i = \f^{(\ell)}\right)\prob\left(\f^{*}_i=\f^{(\ell)}\right)\geq \prob\left(\left.\bhf_i \neq \f^{(j)}\right|\f^{*}_i = \f^{(j)}\right)\prob\left(\f^{*}_i=\f^{(j)}\right). \label{eqn:peExpansion}
\end{align}
}
The denominator term in \eref{eqn:pFDR_oneTest} can be expanded as follows:
\begin{align*}
\prob\left(\bhf_i \neq \f^{(j)} \right) &= \prob\left(\left.\bhf_i \neq \f^{(j)} \right| \f^{*}_i = \f^{(j)}\right)\prob\left(\f^{*}_i = \f^{(j)}\right) \\
&\qquad \qquad+ \prob\left(\left.\bhf_i \neq \f^{(j)} \right| \f^{*}_i \neq \f^{(j)}\right)\prob\left(\f^{*}_i \neq \f^{(j)}\right).
\end{align*}
Observe that $\prob\left(\left.\bhf_i \neq \f^{(j)} \right| \f^{*}_i = \f^{(j)}\right)$ is nonnegative, and
\begin{align*}
&\prob\left(\left.\bhf_i \neq \f^{(j)} \right| \f^{*}_i \neq \f^{(j)}\right) = \prob\left(\left.\bhf_i \in \sD \backslash \f^{(j)}\right|\f^{*}_i\neq \f^{(j)} \right)\geq \prob\left(\left.\bhf_i = \f^{*}_i\right|\f^{*}_i\neq \f^{(j)}\right)\\
&= 1-\prob\left(\left.\bhf_i \neq \f^{*}_i\right|\f^{*}_i \neq \f^{(j)}\right) = 1-\frac{\prob\left(\bhf_i \neq \f^{*}_i,\f^{*}_i \neq \f^{(j)}\right)}{\prob\left(\f^{*}_i \neq \f^{(j)}\right)} \geq 1-\frac{\prob\left(\bhf_i \neq \f^{*}_i\right)}{\prob\left(\f^{*}_i 
\neq \f^{(j)}\right)} \\
&= 1-\frac{\pe_i}{1-p^{(j)}} .
\end{align*}
{
Thus
\begin{align}
\prob\left(\bhf_i \neq \f^{(j)} \right) &\geq \left(1-\frac{\pe_i}{1-p^{(j)}}\right) \left(1-p^{(j)}\right)= 1-p^{(j)} - \pe_i. \label{eqn:pFDRDr_bound}
\end{align}
Substituting \eref{eqn:peExpansion} and \eref{eqn:pFDRDr_bound} in \eref{eqn:pFDR_oneTest},
\begin{align}
\prob\left(\sH_{0i}^{(j)}\left|\by_i \in \Gamma_i^{(j)}\right.\right) &\leq \frac{\pe_i }{1-p^{(j)} - \pe_i} 
\leq \frac{\pe_{\max} }{1-p^{(j)} - \pe_{\max}}. \label{eqn:pFDR_bound}
\end{align}
By substituting \eref{eqn:pFDR_bound} in \eref{eqn: pFDRSubsetEqn} and \eref{eqn:expectation_expn} we have:
\begin{align*}
&\pFDRj\bGammabr \leq \sum_{k=1}^{M} \sum_{S_k} \frac{1}{k}\left(\sum_{\l=1}^{k}\frac{\pe_{\max} }{1-p^{(j)} - \pe_{\max}}\right) \times  \prob \left(\left.\bT \in \Lambda_{S_k} \right| R\bGammabr > 0 \right)\\
&= \frac{\pe_{\max} }{1-p^{(j)} - \pe_{\max}} \sum_{k=1}^{M} \sum_{S_k}\prob \left(\left.\bT \in \Lambda_{S_k} \right| R\bGammabr > 0 \right)\leq \frac{\pe_{\max} }{1-p^{(j)} - \pe_{\max}}
\end{align*}
since $\sum_{k=1}^{M} \sum_{S_k}\prob \left(\left.\bT \in \Lambda_{S_k} \right| R\bGammabr > 0 \right) \leq 1$. The result of Theorem~\ref{thm:pFDR} is obtained by finding an upper bound on the worst-case \pFDR given by
\begin{align*}
&\pFDR_{\max} = \max_{j \in \{1,\ldots,m\}} \pFDRj\bGammabr \\
&\leq \max_{j \in \{1,\ldots,m\}} \frac{\pe_{\max} }{1-p^{(j)} - \pe_{\max}} = \frac{\pe_{\max} }{1-\pmax - \pe_{\max}}
\end{align*}
where $\pmax = \max_{\l \in \{1,\ldots,m\}} p^{(\l)}$.
}

\section{Proof of Theorem~\ref{thm:achievability}}
\label{app:achievability}

The proof is via a random selection technique, similar to random coding arguments common in information theory. Specifically, we will draw a $K \times N$ sensing matrix $\bA$ at random from a particular distribution and then show that, for $\epsilon$, $N$, and $K$ satisfying the conditions of the theorem, the probability that the conclusions of the theorem will fail to hold for this randomly chosen $\bA$ is strictly smaller than unity. This will imply that the conclusions of the theorem must be true for at least one (deterministic) realization of $\bA$.

We begin by specifying all the relevant random variables:
\begin{itemize}
\item $\f^{*}_1,\ldots,\f^{*}_M$ are i.i.d.\ random variables taking values in the dictionary $\sD = \{ \f^{(1)},\ldots,\f^{(m)} \}$ with probabilities $p^{(j)} = \Pr\{ \f^*_i = \f^{(j)} \}, j \in \{1,\ldots,m\}$;
\item $\bn_1,\ldots,\bn_M \stackrel{\text{i.i.d.}}{\sim} \sN(\zeros,\bI)$;
\item $\bG$ is a random $K \times N$ matrix with i.i.d.\ $\sN(0,1)$ entries.
\end{itemize}
We assume that $\{ \f^*_i \}^M_{i=1}$, $\{ \bn_i \}^M_{i=1}$, and $\bG$ are mutually independent, and we will denote by $\PP$ their joint probability distribution. Finally, we let $\bA = \frac{1}{\sqrt{K}}\bG$ and consider the observation model
\begin{align}\label{eq:obs_model}
	\by_i = \alpha_i \bA \f^{*}_i + \bn_i, \qquad i \in \{1,\ldots,M\}
\end{align}
where $\alpha_1,\ldots,\alpha_M > 0$ are the given signal strengths.

We first consider the case when $\alpha_1 = \ldots = \alpha_M = \alpha$. Given $\epsilon$, $N$, and $K$, we define the following two error events:
\begin{align*}
	\sE_1 \deq \left\{ \| \bG \| \ge (1+\epsilon)(\sqrt{K}+\sqrt{N})\right\}, \text{ and }
	\sE_2 \deq \left\{ \bhf_1 \neq \f^{*}_1 \right\},
\end{align*}
where, for each $i \in \{1,\ldots,M\}$, $\bhf_i$ is defined according to \eqref{eqn:fhat}. Note that, since we have assumed that the $\alpha_i$'s are equal and all the pairs $(\f^*_i,\bn_i), i \in\{1,\ldots,M\}$, are i.i.d.,
	\begin{align}\label{eq:p_e_equality}
		\PP(\bhf_i \neq \f^*_i | \bA) = \PP(\sE_2|\bA), \qquad \forall i \in \{1,\ldots,M\}.
	\end{align}
We will now prove that
\begin{align}
	&\PP(\sE_1 \cup \sE_2) \le \frac{1-\pmin}{\pmin} \left(1+\frac{\alpha^2 \dmin^2}{4K\sigma^2}\right)^{-\frac{K}{2}}+ 2\exp\left(-\frac{(K+N)\epsilon^2}{2}\right). \label{eq:p_E}
\end{align}
The union bound gives $\PP(\sE_1 \cup \sE_2) \le \PP(\sE_1) + \PP(\sE_2)$. First, we bound $\PP(\sE_1)$. To do that, we use the following concentration result for Gaussian random matrices \cite{Davidson_Szarek}: for any $t \ge 0$,
\begin{align*}
	\Pr \left\{ \| \bG \| \ge \sqrt{K} + \sqrt{N} + t \right\} \le 2e^{-t^2/2}.
\end{align*}
Letting $t = \epsilon(\sqrt{K}+\sqrt{N})$ and using the fact that $t^2 \ge (K+N)\epsilon^2$, we get
\begin{align}\label{eq:p_E1}
	\PP(\sE_1) \le 2\exp\left(- \frac{(K+N)\epsilon^2}{2}\right).
\end{align}
Next, we bound $\PP(\sE_2)$. To that end, we use the following result, which is a straightforward extension of Theorem~1 in \cite{Nowak_CS_ForSignalClassification} to nonequiprobable dictionary elements:

\begin{lemma}[Compressive classification error] Consider the problem of classifying a signal of interest $\f^{*} \in \sD = \{\f^{(1)},\ldots,\f^{(m)}\}$ to one of $m$ known target classes by making observations of the form $\by = \alpha \bA \f^* + \bn$ where $\bn \sim \Gaussian{\zeros,\sigma^2\bI}$, given the knowledge of the dictionary $\sD$, prior probabilities $p^{(j)}$ for $j \in \{1,\cdots,m\}$, sensing matrix $\bA$, and the noise variance $\sigma^2$. If the entries of $\bA$ are drawn i.i.d. from $\Gaussian{0,1/K}$ independently of $\f^{*}$ and $\bn$, and the estimate $\bhf$ is obtained according to \eref{eqn:fhat}, then 
\begin{align*}
\prob\left(\bhf \neq \f^{*}\right) &\leq \frac{1-\pmin}{\pmin} \left(1+\frac{\alpha^2 \dmin^2}{4K\sigma^2}\right)^{-\frac{K}{2}}
\end{align*} 
where the probability is taken with respect to the distributions underlying $\f^*$, $\bA$, and $\bn$.
\end{lemma}
\noindent Using the above lemma, we have
\begin{align}\label{eq:p_E2}
	\PP(\sE_2) \le \frac{1-\pmin}{\pmin} \left(1+\frac{\alpha^2 \dmin^2}{4K}\right)^{-\frac{K}{2}}.
\end{align}
Combining \eqref{eq:p_E1} and \eqref{eq:p_E2}, we get \eqref{eq:p_E}.

Because of \eqref{eq:positive_prob}, the right-hand side of \eqref{eq:p_E} is less than $1 - \epsilon - p_{\max}$, which is strictly positive by hypothesis. Thus, from the fact that 
\begin{align*}
	\PP(\sE_1 \cup \sE_2) = \expect[ \PP(\sE_1 \cup \sE_2 | \bA)]
\end{align*}
and from \eqref{eq:p_e_equality}, it follows that there exists at least one deterministic choice of the $K \times N$ sensing matrix $\bA^*$, such that:
\begin{subequations}
\begin{align}
\| \bA^* \| &\le (1+\epsilon)\left(1+\sqrt{\frac{N}{K}}\right) \label{eq:A_star_1}\\
\pe_{\max}(\bA^*) &\le \frac{1-\pmin}{\pmin} \left(1+\frac{\alpha^2 \dmin^2}{4K}\right)^{-\frac{K}{2}}
+ 2\exp\left(-\frac{(K+N)\epsilon^2}{2}\right)  \label{eq:A_star_2}
\end{align}
\end{subequations}
where, for a given choice of $\bA$, $\pe_{\max}(\bA)$ denotes the maximum probability of error defined in Theorem~\ref{thm:pFDR}.

Next, from \eqref{eq:A_star_1} and \eqref{eq:weak_background} it follows that $\bA^*$ satisfies the conditions of Theorem~\ref{thm:PhiConstruction}. Finally, we use \eqref{eqn:pfdr_max_anyA} to bound the worst-case \pFDR achievable with $\bA^*$. First of all, we note that the function $U(x) = \frac{x}{1-p_{\max}-x}$ is twice differentiable and convex on the interval $[0,1-p_{\max}]$. Therefore, for any $x \in [0,1-p_{\max}]$ and any $h > 0$ small enough so that $x+h \in [0,1-p_{\max}]$, we have
\begin{align}
	U(x+h) &\le U(x) + U'(x+h)h= U(x) + \frac{(1-p_{\max})h}{(1-p_{\max}-x-h)^2}. \label{eq:first_order_approx}
\end{align}
Let us choose
\begin{align*}
x =  \frac{1-\pmin}{\pmin} \left(1+\frac{\alpha^2 \dmin^2}{4K}\right)^{-\frac{K}{2}} \text{ and }
 h = 2\exp\left(-\frac{(K+N)\epsilon^2}{2}\right).
\end{align*} 
Then from \eqref{eq:positive_prob} we have $x+h \le 1-\epsilon-p_{\max} < 1 - p_{\max}$, and from \eqref{eqn:boundonK} we have $x+h \ge 0$. Hence, using \eqref{eq:first_order_approx} and simplifying, we obtain the bound
\begin{align*}
	&\pFDR_{\max}(\bA^*) \le \frac{1}{\pmin}\left(\frac{1-\pmax}{1-\pmin}\left(1+\frac{\alpha^2 \dmin^2}{4K}\right)^{\frac{K}{2}}-\frac{1}{\pmin}\right)^{-1} + \\
	&\qquad\frac{2(1-\pmax)}{\epsilon^2}\exp\left(-\frac{(K+N)\epsilon^2}{2}\right).
\end{align*}
This proves the theorem for the case $\alpha_1 = \ldots = \alpha_M = \alpha$.

To handle the case when the $\alpha_i$'s are distinct, we simply let
$$
i^* \deq \argmin_{i \in \{1,\ldots,M\}} \alpha_i
$$
and replace the definition of the error event $\sE_2$ with $\sE'_2 = \{ \bhf_{i^*} \neq \f^{*}_{i^*}\}$. Then the same argument goes through, except that instead of \eqref{eq:p_e_equality} we use the bound
\begin{align*}
	\PP(\bhf_i \neq \f^*_i|\bA) \le \PP(\bhf_{i^*} \neq \f^*_{i^*}|\bA) = \PP(\sE'_2|\bA), \,\, \forall i \neq i^*
\end{align*}
which follows from the following argument. First of all, we can replace the observation model with the equivalent model
\begin{align*}
	\bty_i = \bA \f^*_i + \btn_i, \qquad i \in \{1,\ldots,M\}
\end{align*}
where $\btn_i = \frac{1}{\alpha_i}\bn_i \sim \sN(\zeros,\frac{1}{\alpha^2_i}\bI)$. Secondly, from the fact that $\alpha_i \ge \alpha_{i^*} \equiv \alpha_{\min}$ for any $i \neq i^*$ it follows that $\btn_{i^*}$ is equal in distribution to $\btn_i + \btn'_i$, where $\btn'_i \sim \sN\left(\zeros, \left(\frac{1}{\alpha_i^2} - \frac{1}{\alpha^2_{\min}}\right)\bI\right)$ is independent of $\btn_i$. This implies that the $i^*$th observation is the noisiest, and the corresponding MAP estimate $\bhf_{i^*}$ has the largest probability of error.

\section{Proof of Theorem~\ref{theorem:AnomDet}}
\label{app:AnomDet}
We first prove this theorem assuming that $\{\alpha_i\}$ are known and
later extend to the case where $\{\halpha_i\}$ are estimated from the
observations. Let $\btf_i = \argmin_{\f \in \sD}\|\f^{*}_i-\f\|$. The
p-value expression in \eref{eqn:pvalue} can be expanded as follows:
\begin{align}
&p_i = \prob\left(\left.{\td_i} \geq d_i \right| \sH_{0i}\right)= \prob\left(\left. \min_{\f \in \sD} \|\alpha_i\bA(\f^{*}_i-\f)+\bn\|\geq d_i \right|\sH_{0i}\right) \nonumber\\
&\leq \prob\left(\left. \|\alpha_i\bA(\f^{*}_i-\btf_i)+\bn\|\geq d_i \right|\sH_{0i}\right) = \prob\left(\left. \|\alpha_i\bA(\f^{*}_i-\btf_i)+\bn\|^2\geq d_i^2 \right|\sH_{0i}\right)
\label{eqn:pvalueExpn_expanded}.
\end{align}
Note that $\|\alpha_i\bA(\f^{*}_i-\btf_i)+\bn\|^2$ is a noncentral $\chi^2$ random variable with $K$ degrees of freedom and a noncentrality parameter $\nu_{i} = \|\alpha_i\bA\left(\f_i^{*}-\btf_i\right)\|^2$. Thus \eref{eqn:pvalueExpn_expanded} can be written in terms of a noncentral $\chi^2$ CDF $\sF\left(d_i^2; K,\nu_{i}\right)$ with parameter $d_i^2$. The upper and lower bounds on $\nu_{i}$ can be obtained using the properties of the projection matrix $\bA$. 
Applying \eref{eqn:PrePros_distPres}, we see that
\begin{align*}
\alpha_i^2(1-\epsilon)^2\|\f^{*}_i-\btf_i\|^2 \leq \nu_{i}\leq \alpha_i^2(1+\epsilon)^2\|\f^{*}_i-\btf_i\|^2
\end{align*}
with high probability. 
Thus, 
\begin{align}
&p_i \leq 1-\prob\left(\left. \|\alpha_i\bA(\f^{*}_i-\btf_i)+\bn\|^2\leq d_i^2 \right|\sH_{0i}\right) = 1-\sF\left(d_i^2; K,\nu_{i}\right) \label{eqn:pvalueUB} \\
&\leq 1-\sF\left(d_i^2; K,\alpha_i^2(1+\epsilon)^2\|\f^{*}_i-\btf_i\|^2\right) \leq 1-\sF\left(d_i^2; K,\alpha_i^2(1+\epsilon)^2\tau^2\right) \nonumber
\end{align}
since $\|\f^{*}_i-\f\| \leq \tau$ for all $\f \in \sD$ under $\sH_{0i}$. 

When $\{\alpha_i\}$ are estimated from the observations such that $\{\halpha_i\}$ satisfy \eref{eqn:estAccuracy}, we can write the p-value expression in \eref{eqn:pvalueUB} as follows: 
\begin{align}
p_i &\leq 1-\sF\left(d_i^2; K,\left\|\bA\left(\alpha_i\f^{*}_i-\halpha_i\btf_i\right)\right\|^2\right) \nonumber \\
&\leq 1-\sF\left(d_i^2; K,(1+\epsilon)^2\halpha_i^2\left\|\frac{\alpha_i}{\halpha_i}\f^{*}_i-\btf_i\right\|^2\right) \label{eqn:distPres_2}
\end{align}
where{\eref{eqn:distPres_2}} is due to the distance preservation property of $\bA${given in \eref{eqn:PrePros_distPres}}. Observe that $\left\|\frac{\alpha_i}{\halpha_i}\f^{*}_i-\btf_i\right\|^2$ can be upper bounded as shown below: 
\begin{align*}
\left\|\frac{\alpha_i}{\halpha_i}\f^{*}_i-\btf_i\right\|^2 &= \left\|\left(\frac{\alpha_i}{\halpha_i}-1\right)\f^{*}_i+\f^{*}_i-\btf_i\right\|^2 \leq \left(\left\|\left(\frac{\alpha_i}{\halpha_i}-1\right)\f^{*}_i\right\|+\big\|\f^{*}_i-\btf_i\big\|\right)^2 \\
&= \left(\left|\frac{\alpha_i}{\halpha_i}-1\right|+\big\|\f^{*}_i-\btf_i\big\|\right)^2 \leq \left( \zeta+\big\|\f^{*}_i-\btf_i\big\|\right)^2 \end{align*}
where third-to-last equation is due to the triangle inequality, second-to-last equation comes from the assumption that $\left\|\f^{*}_i\right\|=1$, and the last inequality is due to \eref{eqn:estAccuracy}. By applying this result to \eref{eqn:distPres_2} and exploiting the fact that $\|\f^{*}_i-\f\|\leq\tau$ under $\sH_{0i}$ for some $\f \in \sD$, we have
\begin{align*}
p_i&\leq 1-\sF\left(d_i^2; K,(1+\epsilon)^2\halpha_i^2\left(\zeta+\big\|\f^{*}_i-\btf_i\big\|\right)^2\right) \leq 1-\sF\left(d_i^2; K,(1+\epsilon)^2\halpha_i^2\left(\zeta+\tau\right)^2\right) .
\end{align*}

\end{appendix}
\bibliographystyle{plain}
\bibliography{TIP2010.bib}

\end{document}

%% file: myMacros_arxiv.tex
\usepackage{amsmath,amssymb,latexsym}

\usepackage{stmaryrd}
\usepackage{amsfonts}
\usepackage{bm,xspace}
\usepackage{bibunits}
\usepackage{bbm}
\usepackage[pdftex]{graphicx}

\usepackage[usenames]{color}

\newcommand{\setbr}[1]{{\left\{{#1}\right\}}}
\newcommand{\paranbr}[1]{{\left({#1}\right)}}
\newcommand{\squarebr}[1]{{\left[{#1}\right]}}
\newcommand{\Gaussian}[1]{{\sN\left({#1}\right)}}
\newcommand{\abs}[1]{{\left|{#1}\right|}}
\newcommand{\norm}[1]{{\left\|{#1}\right\|}}

 \newcommand{\xmath}[1]{{\ensuremath{#1}\xspace}}
 \newcommand{\bmath}[1]{\xmath{\bm{#1}}}

 \def \sE{{\mathcal E}}

 \def \tGamma{\widetilde{\Gamma}}

 \def \sD{\mathcal{D}}

 \def \l{\ell}

 \def \f{\bmath{f}}
 \def \bx{\bmath{x}}

 \def \bg{\bmath{g}}
 \def \bG{\bmath{G}}

 \def \bhf{\bmath{\widehat{f}}}
 \def \btf{\bmath{\widetilde{f}}}
 
 \def \by{\bmath{y}}
 \def \bn{\bmath{n}}
 \def \bz{\bmath{z}}

 \def \ie{{\em i.e.},\ }

\def \th{^{\rm th}}

\def \sN{{\mathcal N}}

\def\ith{{\xmath{i^{{\rm th}}}}\xspace}
\def\jth{{\xmath{j^{{\rm th}}}}\xspace}

 \newcommand{\bmu}{\bmath{\mu}}

\newcommand{\bPhi}{\bmath{\Phi}\xspace}

\newcommand{\bSigma}{\bmath{\Sigma}\xspace}

\newcommand{\btSigma}{\widetilde{\bSigma}\xspace}

 \def\argmin{\mathop{\rm arg\,min}}
 \def\argmax{\mathop{\rm arg\,max}}
 \newcommand{\reals}{{\xmath{\mathbb{R}}}}
 
 \newcommand{\expect}{{\xmath{\mathbb{E}}}}
 
\newcommand{\prob}{{\xmath{\mathbb{P}}}}
 \newcommand{\ind}[1]{{\xmath{\mathbb{I}_{\left\{#1\right\}}}}}
 \newcommand{\eref}[1]{(\ref{#1})}

 \newcommand{\Order}{\mathcal{O}}

\newcommand{\ty}{\widetilde{y}}
\newcommand{\bty}{\bmath{\ty}}

\newcommand{\sF}{\mathcal{F}}
\newcommand{\sH}{\mathcal{H}}
\newcommand{\sM}{\mathcal{M}}

\newcommand{\bI}{\bmath{I}}

\newcommand{\bw}{\bmath{w}}
\newcommand{\bA}{\bmath{A}}

\newcommand{\zeros}{\bmath{0}}
\newcommand{\identity}{{\xmath{\mathbb{I}}}}

\newcommand{\bB}{\bmath{B}}
\newcommand{\bb}{\bmath{b}}

\newcommand{\bT}{\bmath{T}}
\newcommand{\halpha}{\widehat{\alpha}}

\newcommand{\pFDR}{{\rm{pFDR}}\xspace}
\newcommand{\FDR}{{\rm{FDR}}\xspace}
\newcommand{\FNR}{{\rm{FNR}}\xspace}

\newcommand{\setGamma}{\left(\left\{\Gamma_i\right\}\right)\xspace}

\newcommand{\pe}{({\rm{P_e}})}

\newcommand{\pmin}{p_{\min}}
\newcommand{\pmax}{p_{\max}}
\newcommand{\dmin}{d_{\min}}
\newcommand{\pFDRj}{\pFDR^{(j)}}
\newcommand{\bGammabr}{\left(\bm{\Gamma}\right)\xspace}

\newcommand{\td}{\widetilde{d}}
\newcommand{\sU}{\mathcal{U}}
\newcommand{\btg}{\widetilde{\bg}}
\newcommand{\tg}{\widetilde{g}}
\newcommand{\btb}{\widetilde{\bb}}
\newcommand{\fstar}{\f^*}
\newcommand{\btn}{\widetilde{\bn}}

\def\deq{\stackrel{\scriptscriptstyle\triangle}{=}}

\newcommand{\squishlist}{
   \begin{list}{$\bullet$}
    { \setlength{\itemsep}{0pt}      \setlength{\parsep}{0pt}
      \setlength{\topsep}{0pt}       \setlength{\partopsep}{0pt}
      \setlength{\leftmargin}{1.5em} \setlength{\labelwidth}{1em}
      \setlength{\labelsep}{0.5em} } }

\newcommand{\squishlisttwo}{
   \begin{list}{$\bullet$}
    { \setlength{\itemsep}{0pt}    \setlength{\parsep}{0pt}
      \setlength{\topsep}{0pt}     \setlength{\partopsep}{0pt}
      \setlength{\leftmargin}{2em} \setlength{\labelwidth}{1.5em}
      \setlength{\labelsep}{0.5em} } }

\newcommand{\squishend}{
    \end{list}  }